\documentclass[11pt]{article}
\usepackage{amssymb}
\usepackage{babel}[english]
\usepackage{amsmath,amsfonts,amsthm,amssymb}
\usepackage[dvips]{graphics}
\usepackage{epsfig}
\usepackage{indentfirst}
\usepackage{color}
\usepackage{amsmath}
\usepackage{relsize}
\usepackage{comment}
\usepackage{hyperref}
\usepackage{enumerate}
\allowdisplaybreaks

\newtheoremstyle{theorem}
  {10pt}          
  {10pt}  
  {\sl}  
 {}
  {\bf}  
  {. }    
  { }    
  {}     
\theoremstyle{theorem}

\newtheorem{theorem}{Theorem}[section]
\newtheorem{corollary}{Corollary}[section]
\newtheorem{definition}{Definition}[section]
\newtheorem{lemma}{Lemma}[section]
\newtheorem{proposition}{Proposition}[section]
\newtheorem{remark}{Remark}[section]

\numberwithin{equation}{section}

\newtheoremstyle{defi}
{10pt}  
{10pt}  
{\rm}   
{}      
{\bf}   
{. }    
{ }     
{}      
\theoremstyle{defi}

\def\dd{\mathrm d}
\def\eps{\varepsilon}

\def\Q{\mathbf Q}

\renewcommand{\div}{\mathrm{div}\,}
\newcommand{\R}{\mathbb{R}}

\newcommand{\seps}{\sigma_{\eps}}

\newcommand{\reps}{\rho_{\eps}}

\newcommand{\meps}{\mathbf{m}_{\eps}}
\newcommand{\tmeps}{\widetilde{\mathbf{m}}_{\eps}}
\newcommand{\ueps}{u_{\eps}}
\newcommand{\phie}{\phi_{\eps}}
\newcommand{\Feps}{F_{\eps}}
\newcommand{\tFeps}{\tilde{\Feps}}

\newcommand{\Heps}{H_{\eps}}
\newcommand{\eith}{e^{it\Heps}}

\newcommand{\loc}{\rm loc}
\newcommand{\ue}{\mathbf{u}^E}
\newcommand{\supp}{\operatorname{supp}}

\definecolor{mypink}{RGB}{219, 48, 122}
\definecolor{myblue}{RGB}{0, 0, 122}

\begin{document}

\date\today

\title{Inviscid incompressible limit for capillary fluids with density dependent viscosity}

\author{Matteo Caggio$^{1}$ \ \ Donatella Donatelli$^{2}$ \ \ 
Lars Eric Hientzsch$^{3}$
\\
{\small  1. Institute of Mathematics of the Academy of Sciences of the Czech Republic,} \\
{\small \v Zitn\' a 25, 11567, Praha 1, Czech Republic}\\
{\small  2. Department of Information Engineering, Computer Science and Mathematics,}\\
{\small University of L'Aquila}\\
{\small via Vetoio, Coppito - 67100 L' Aquila, Italy}\\
{\small 3. Department of Mathematics, Karlsruhe Institute of Technology,}\\
\small{Englerstraße 2, 76131 Karlsruhe, Germany}\\  
{\small caggio@math.cas.cz}\\
{\small donatella.donatelli@univaq.it}\\
{\small lars.hientzsch@kit.edu}
}
\date{\today}

\maketitle

\begin{abstract}
    \noindent
	The asymptotic limit of the  Navier-Stokes-Korteweg system for barotropic capillary fluids with density dependent viscosities in the low Mach number and vanishing viscosity regime is established in $\mathbb{R}^d$, with  $d=2,3$. In the relative energy framework, we prove the convergence of weak solutions of the Navier-Stokes-Korteweg system to the strong solution of the incompressible Euler system. The convergence is obtained through the use of suitable dispersive estimates for an acoustic system altered by the presence of the Korteweg tensor.
\end{abstract}
\smallskip

\noindent\textbf{Key words}: barotropic compressible fluids, density dependent viscosity, Navier-Stokes-Korteweg model, capillary fluids, inviscid incompressible limit.

\tableofcontents{}

\section{Introduction} \label{intr}

The aim of this paper is to characterize the asymptotic limit in the low Mach and vanishing viscosity limit regime of the following (re-scaled) compressible Navier-Stokes-Korteweg system in $(0,T) \times \R^3$:
\begin{gather}
\begin{cases}
\partial_{t}\varrho +\textrm{div}(\varrho \mathbf{u})=0,\\\noalign{\vskip0.5mm}
\label{mom}
\partial_{t} (\varrho \mathbf{u})+\textrm{div}(\varrho\mathbf{u}\otimes\mathbf{u})
 +  \frac{\nabla p(\varrho)}{\varepsilon^2}
 -2\nu \textrm{div} (\varrho \mathbb{D}(\mathbf{u}))
 - 2\kappa^2 \varrho\nabla\Delta\varrho=0,
 \end{cases}
\tag{NSK}
\end{gather}
complemented with far-field behavior 
\begin{equation}\label{farfield}
    \varrho\rightarrow 1, \quad \sqrt{\rho}\mathbf{u}\rightarrow 0,  \quad \text{as} \,\, |x|\rightarrow \infty,
\end{equation}
and the initial conditions
\begin{equation} 
\label{ic} \varrho (0,\cdot) =
\varrho_{0}, \qquad \varrho \mathbf{u} (0,\cdot) = \varrho_{0}\mathbf{u}_{0}.
\end{equation}

\noindent
The unknown variables $\varrho=\varrho\left(t,x\right)$,
$\mathbf{u}=\mathbf{u}\left(t,x\right)$ and
$p=p(\varrho)$ represent the mass density,
the velocity vector and the pressure of the fluid, respectively. This
last is given by a standard power law type
\begin{equation} 
\label{press} 
p(\varrho) = \varrho^\gamma, \quad\,\,\ \gamma > 1,
\end{equation}
while 
$$\mathbb{D}(\mathbf{u})=\frac12(\nabla \mathbf{u} +
\nabla^{\top} \mathbf{u})$$ is the symmetric part of $\nabla \mathbf{u}$. The Korteweg tensor $2\kappa^2 \varrho\nabla\Delta\varrho$ can be recast in the form $\varrho\nabla\Delta\varrho=\div \mathbb{K}$ where  
\begin{equation} \label{K}
\mathbb{K} =  
\left(\varrho  \text{div}(\nabla \varrho) + \frac{1}{2} |\nabla \varrho |^2
\right)\mathbb{I}
- (\nabla \varrho \otimes \nabla \varrho).
\end{equation}
Fluids for which the Korteweg tensor has the form as in (\ref{K}) are usually termed capillary fluids (see e.g. \cite{PaoloStefano} and references therein). 
The parameter $\varepsilon$ represents the Mach number while $\nu$ and $\kappa$ are the viscosity and capillary coefficients, respectively. The associated energy is given by 
\begin{equation}\label{eq:energy}
    E(\varrho,\mathbf{u})(\tau) = \int_{\R^3}
    \left[\frac{\varrho|\mathbf{u}|^2}{2}+ H(\varrho) +\kappa^2 |\nabla \varrho|^2\right](\tau) dx,
\end{equation}
\begin{equation} \label{H-rho}
    H(\varrho)=\frac{1}{\gamma-1}(\varrho^{\gamma}-1-\gamma(\varrho-1)),
\end{equation}
where $H(\varrho)$ takes the far-field behavior into account. 

\subsection{Singular limit, scope of the present analysis} \label{sing}

Formally, in the limit of $\varepsilon\to 0$ and $\nu\to 0$, we obtain the incompressible Euler system in the whole space

\begin{gather}\label{EULER}
\begin{cases}
    \ \ \text{div}_{x}\mathbf{u}^E =0, \\ \partial_t \mathbf{u}^E + \mathbf{u}^E \cdot \nabla \mathbf{u}^E + \nabla \Pi =0,
\end{cases}
\tag{E}
\end{gather}
with initial conditions
\begin{equation} \label{EULER-ic}
    \mathbf{u}^E (0,\cdot) = \mathbf{u}_{0}^E.
\end{equation}
More precisely, in the present analysis we are interested to prove the convergence of the weak solution of \eqref{mom} to the strong solution of \eqref{EULER} in the limit of $\varepsilon\to 0$ and $\nu\to 0$, namely in the low-Mach (incompressible) and vanishing viscosity limit. The key difficulty of our analysis consists in tackling the singular limit in the presence of a highly non-linear capillarity tensor, a density-dependent viscosity and disturbing effects caused by fast oscillating acoustic waves. The occurrence of the latter is due to the choice of the so-called ill-prepared initial data (see Section \ref{MR}).

The convergence will be obtained within the relative energy inequality framework, we refer for instance to \cite{FN} and references therein for a comprehensive overview of the method. To the best of our knowledge, the present is the first result addressing the inviscid incompressible limit for capillary fluids and could be seen as a continuation of a previous analysis in which a weak-strong uniqueness result (for a fixed viscosity coefficient $\nu > 0$), together with the high-Mach number limit,  has been recently obtained; see \cite{caggio-donatelli-2}.

\bigskip
\noindent
To achieve our result, we thus need the following:
\begin{itemize}
    \item A suitable form of the relative energy inequality that captures structural properties of compressible flow including the capillarity and density-dependent viscosity tensors, in particular, in the vanishing viscosity limit.
    \item Appropriate dispersive estimates for the decay of the acoustic waves in the low-Mach number limit that take the capillarity effects into account.
\end{itemize}

Regarding the first point, in \cite{caggio-donatelli-2} the authors derived a relative energy inequality based on an ``augmented" version of the system (\ref{mom}) in the same spirit of e.g. \cite{bresch-3}, \cite{bresch-1}, \cite{bresch-2}, where an equation for a velocity of the type $2\nu\nabla \log(\varrho)$ is properly introduced.
The use of the ``augmented" system is required
due to the fact that an $H^1$ bound for the velocity is no longer available because of the density dependent viscosity. Consequently, standard application of the Korn’s inequality to handle, in particular, the viscous terms 
in the usual weak-strong uniqueness framework is not possible.
The authors in \cite{BC} used this ``augmented" version to study, for the first time, the vanishing viscosity limit for a barotropic fluid with density dependent viscosity. The result has been obtained for a fluid in smooth bounded domain and without the presence of the Korteweg term. However, a recent analysis (see \cite{BCD}) shows that it is possible to completely avoid the ``augmented" version of the system \ref{mom} using a suitable form of the relative energy inequality. 
Consequently, being the proper framework to perform the vanishing viscosity limit, our relative energy inequality will be consistent with the one considered in \cite{BCD}.

On the other hand, the analysis of the incompressible limit requires a suitable control of the acoustic waves which presence is due to (potential) density fluctuations at finite Mach number. Being the problem \eqref{mom} posed on $\R^3$, suitable dispersive estimates yield the decay of the acoustic waves in the low-Mach number limit. Here, we take the presence of the nonlinear capillarity tensor $2\kappa^2\rho\nabla\Delta \rho$ depending on third order derivatives of the density into account. Upon linearizing system \eqref{mom} around the constant state $(\rho=1, u=0)$, it is possible to obtain a fourth-order acoustic wave equation for the density fluctuations $\seps=\eps^{-1}(\varrho_\varepsilon-1)$, namely
\begin{equation}\label{eq:densityfluct-intro}
    \partial_{tt}^2 \seps-\frac{1}{\eps^2}\Delta(1-2\eps^2\kappa^2\Delta) \seps=0.
\end{equation}
The presence of the Korteweg tensor modifies the dispersion relation governing the propagation of the acoustic waves with respect to the usual dispersive estimates obtained for the classical wave equation, namely in the absence of the capillarity tensor $(\kappa=0)$ and is reminiscent to the accoustic oscillations occuring in quantum fluids \cite{lars1, lars3}.
Consequently, the adaptation of the techniques developed in \cite{lars1, lars3} allows us to derive suitable estimates for the present analysis (see Section \ref{AWe} below).

\subsection{Overview of previous results}
Other results in the context of Korteweg fluids concern the low-Mach number limit for quantum fluids (see \cite{lars1}, \cite{lars3}, \cite{lars2}) in which the authors study a weak-weak convergence towards the weak solution of the quantum incompressible Navier-Stokes system. 

A similar singular limit analysis for the Korteweg type fluids has also been performed in the case of the quasi neutral limit; see \cite{bresch05},\cite{donatella-piero}.
In the framework of strong solutions we can refer to \cite{FuLi}, \cite{HaLiRo}, \cite{LiYo},  \cite{JuXu}, \cite{ShaLi}, \cite{LY15}.

\subsection{Organization of the paper}
The manuscript is organized as follows: in Section \ref{background} we discuss the existence of weak solutions for the Navier-Stokes-Korteweg system (\ref{mom})
together with the existence of the strong solution of the Euler system (\ref{EULER}). Then, we present our main result. Section \ref{ID-UE} is devoted to the uniform bounds. In Section \ref{AWe} we introduce the acoustic system related to the Navier-Stokes-Korteweg system (\ref{mom}) and the dispersive estimates. In Section \ref{REI} we derive the relative energy inequality suitable for our analysis. In Section \ref{IIL}, we study the incompressible and vanishing viscosity limit in the whole space $\mathbb{R}^3$. Section \ref{2D-c} is devoted to some comments and observation that concern the incompressible and vanishing viscosity limit in the 2D-case.   

\subsection{Notation} \label{not}

We denote by $C_c^\infty([0,T)\times\R^3;\R^d)$ the space of periodic smooth functions with values in $\R^d$ with compact support in $[0,T)\times\R^3$ and by $L^p(\R^3)$ the standard Lebesgue spaces. The Sobolev spaces of functions with $s$ distributional derivative in $L^p(\R^3)$ are $W^{s,p}(\R^3)$. In the case $p=2$, $W^{s,p}(\R^3) = H^s(\R^3)$. The Bochner spaces for time-dependent functions with values in Banach spaces $X$ are denoted by $L^p(0,T;X)$ and $W^{k,p}(0,T;X)$. The space $C(0,T;X_w)$ is the space of continuous functions endowed with the weak topology. By $\mathbf{Q}$ and $\mathbf{P}$ we denote the Helmholtz-Leray projectors on irrotational and divergence-free vector fields, respectively:
$$
\mathbf{Q} = \nabla \Delta^{-1}\text{div}, \ \ \ \ \mathbf{P}=\mathbf{I}-\mathbf{Q}.
$$

\section{Background and Main Result}
\label{background}

We introduce the concepts of weak solutions for the Navier Stokes Korteweg system \eqref{mom} and of strong solutions for the target Euler system \eqref{EULER}. For both systems solutions are considered in $\R^3$ with given farfield condition \eqref{farfield}.

\subsection{Existence of weak solutions for capillary fluids} \label{ws-ex}
\noindent
\begin{definition} \label{weak-def}
A triple $(\varrho,\mathbf{u},\mathcal{T})$, with $\varrho \geq 0$, is said to be a weak solution to (\ref{mom}) with initial data (\ref{ic}) if the following conditions are satisfied:

\begin{enumerate}[(i)]
    \item Integrability conditions:
\begin{equation*}
    \varrho\in L^\infty(0,T;H_{\loc}^1(\R^3))\cap
    L^2(0,T;H_{\loc}^2(\R^3)), \ \ \ \ 
    \sqrt{\varrho} \mathbf{u} \in L^\infty(0,T;L_{\loc}^2(\R^3)),
\end{equation*}
\begin{equation*}
    \varrho^{\frac{\gamma}{2}}\in
    L^\infty(0,T;L_{\loc}^2(\R^3))\cap
    L^2(0,T;H_{\loc}^1(\R^3)), \ \ \ \ 
    \nabla \sqrt{\varrho} \in
    L^\infty(0,T;L_{\loc}^2(\R^3)),
\end{equation*}
\begin{equation*}
    \mathcal{T} \in 
    L^2(0,T;L_{\loc}^2(\R^3)), \ \ \ \ 
    \varrho \mathbf{u} \in C([0,T); L_{\loc,w}^{\frac{3}{2}}(\R^3)).
\end{equation*}

\item Continuity equation:
for any $\varphi \in C_c^\infty([0,T]\times\R^3;\R)$ and for all $\tau \in [0,T]$,
\begin{equation} \label{cont-wf}
- \int_0^{\tau} 
\int_{\R^3}
\left(
\varrho \partial_t \varphi + \sqrt{\varrho}\sqrt{\varrho} \mathbf{u} \cdot \nabla \varphi 
\right)
dx dt = \int_{\R^3}
\varrho (0,\cdot) \varphi (0,\cdot) dx - \int_{\R^3}
\varrho (\tau,\cdot) \varphi (\tau,\cdot) dx.
\end{equation}

\item Momentum equation:
for any fixed $l=1,2,3$, $\phi \in C_c^\infty([0,T]\times\R^3;\R^3)$
and for all $\tau \in [0,T]$,

\begin{equation*}  
    - \int_0^{\tau} \int_{\R^3}
    \sqrt{\varrho}\sqrt{\varrho}\mathbf{u} 
    \cdot
    \partial_t \phi \; dx dt
    - \int_0^{\tau} \int_{\R^3}
    (\sqrt{\varrho}\mathbf{u} \otimes \sqrt{\varrho} \mathbf{u}) :\nabla \phi \; dx dt
\end{equation*}
\begin{equation*}
    + 2\nu \int_0^{\tau} \int_{\R^3} 
    \sqrt{\varrho} \mathcal{S}(\mathbf{u}) :
    \nabla \phi \; dx dt
    - \frac{1}{\varepsilon^2}\int_0^{\tau} \int_{\R^3}
    p(\varrho) \text{div}\phi \; dx dt
\end{equation*}
\begin{equation*} 
    -2\kappa^2 \int_0^{\tau} \int_{\R^3} 
    \left( 
    \nabla\varrho \cdot \nabla( \varrho \textrm{div}\phi) 
    - \frac{1}{2} |\nabla \varrho|^2 \textrm{div}\phi
    + \nabla \varrho \otimes \nabla \varrho : \nabla \phi
    \right)dxdt
\end{equation*}
\begin{equation} \label{mom-wf}
    = 
    \int_{\R^3} (\varrho \mathbf{u} \cdot \phi)(0,\cdot) dx
    - \int_{\R^3} (\varrho \mathbf{u} \cdot \phi)(\tau,\cdot) dx.
\end{equation}

\item Dissipation:
for any $\xi \in C_c^\infty([0,T]\times\R^3;\R)$ and for all $\tau \in [0,T]$,

\begin{equation} \label{diss}
    \int_0^{\tau} \int_{\R^3}
    (\sqrt{\varrho} \mathcal{T}) \xi \; dxdt
    =
    - \int_0^{\tau} \int_{\R^3}
    \varrho \mathbf{u} \cdot \nabla \xi \; dxdt
    - \int_0^{\tau} \int_{\R^3} 
    2 (\sqrt{\varrho} \mathbf{u} \otimes \nabla \sqrt{\varrho}) \xi \; dxdt.
\end{equation}

\bigskip
\item Energy inequality:
for almost all $\tau \in [0,T]$, the energy as defined in \eqref{eq:energy} satisfies the inequality 
\begin{equation} \label{eq:EI}
    E(\varrho,\mathbf{u})(\tau)+ 2\nu \int_0^{\tau} \int_{\R^3}
    |\mathcal{S}(\mathbf{u})(\tau,x)|^2 dxdt\leq E(\varrho^0,\mathbf{u}^0),
\end{equation}
where $\mathcal{S}$ denotes the symmetric part of the tensor $\mathcal{T}$.
\item BD-entropy inequality: there exists $C>0$ such that for almost all $\tau \in [0,T]$ the Bresch-Desjardins entropy
inequality holds 
\begin{equation*}
\mathcal{B}(\varrho, \mathbf{u})(\tau) +2\nu\int_{0}^\tau\int_{\R^3}|\mathcal{A}(\mathbf{u})(\tau,x)|^2 dxdt
\end{equation*}
\begin{equation} \label{BD}
+\frac{8\nu}{\gamma^2}\int_{0}^\tau\int_{\R^3}\left|\nabla \varrho^{\frac{\gamma}{2}}(\tau,x)\right|^2 dxdt+4\nu\kappa^2\int_{0}^\tau\int_{\R^3}|\Delta \varrho(\tau,x)|^2 dxdt \leq C \mathcal{B}(\varrho^0, \mathbf{u}^0),
\end{equation}
with
\begin{equation} \label{eq:BDI}
    \mathcal{B}(\varrho,\mathbf{u})(\tau)=\int_{\R^3}\left[\frac{1}{2}|\sqrt{\varrho}u+2\nu\nabla \sqrt{\varrho}|^2+\kappa^2|\nabla{\varrho}|^2+H(\varrho)\right](\tau) dx.
\end{equation}
and $\mathcal{A}$ denoting the anti-symmetric part of $\mathcal{T}$.
\end{enumerate}
\end{definition}

Here, $\mathcal{S}(\mathbf{u})$ and $\mathcal{A}(\mathbf{u})$ are the symmetric and anti-symmetric part of the tensor $\mathcal{T}(\mathbf{u})$, respectively, defined by
\begin{equation} \label{TT}
    \sqrt{\varrho}\mathcal{T}(\mathbf{u})
    = \nabla (\varrho \mathbf{u}) - 2  (\sqrt{\varrho} \mathbf{u}) \otimes \nabla \sqrt{\varrho}
    \ \ \ \ 
    \text{in}
    \ \ \ \ 
    \mathcal{D}'\left([0,T]\times \R^3\right).
\end{equation}

\begin{remark} \label{smooth-E}
For smooth solutions, the energy inequality (in fact an equality) for the system (\ref{mom}) reads 
\begin{equation} \label{ee-sm}
    E(\varrho,\mathbf{u})(\tau) + 2\nu\int_0^{\tau} \int_{\R^3} \varrho|\mathbb{D}(\mathbf{u})|^2 dxdt \leq E(\varrho^0,\mathbf{u}^0).
\end{equation}
For weak solutions, we only require the energy inequality \eqref{eq:EI} to hold where the dissipation \eqref{diss} is formulated in distributional sense. Weak solutions are commonly constructed through approximation procedures and at present it appears to be unclear whether arbitrary finite energy weak solutions satisfy the inequality (\ref{ee-sm}). Indeed, a lack of compactness rules out to conclude that the weak-limit in $L_{t,x}^2$ of $\sqrt{\varrho_n}\mathbb{D}(\mathbf{u}_n)$ is given by $\sqrt{\varrho}\mathbb{D}(\mathbf{u})$, see \cite[Remark 2.2]{PaoloStefano}. In general, the approximation procedure only yields the information
\begin{equation*}
    \sqrt{\varrho_n}\mathbb{D}(\mathbf{u}_n) \rightharpoonup \mathcal{S} \ \ \text{in} \ \ L_{t,x}^2.
\end{equation*}
The viscous term ${\varrho} \mathbb{D}(\mathbf{u})$ has hence to be understood as $\sqrt{\varrho}\mathcal{S}$. Similar arguments hold for the BD-entropy inequality.
\end{remark}

We postulate the following existence result of global finite energy weak solutions. While to the best of the authors' knowledge this result is not available in the literature for \eqref{mom} considered on $\R^3$ with far-field \eqref{farfield}, the respective result for the problem posed on $\mathbb{T}^3$ is proven in \cite{PaoloStefano}. It is conceivable that global existence of weak solutions to \eqref{mom} on $\R^3$ with far-field \eqref{farfield} can be achieved by combining the periodic result \cite{PaoloStefano} and the invading domains approach developed in \cite{LarsPaoloStefano} in the case of quantum fluids in order to prove the following:

\begin{theorem} \label{main-ex}
Given initial data $(\varrho^0, \mathbf{u}^0)$ of finite energy and BD-entropy, then, there exists at least a global weak solution $(\varrho, \mathbf{u}, \mathcal{T})$ of \eqref{mom} in the sense of Definition \ref{weak-def}. 
\end{theorem}

\subsection{Local well-posedness for the incompressible Euler equations}

We recall here the following classical result (see \cite{Kato}, \cite{KL84}  for example) for the target Euler system (\ref{EULER}):

\begin{theorem}\label{3dEuler-ex}
Given $\mathbf{u}_0^E\in W^{3,2}(\R^3)$ with ${\rm div}\mathbf{u}_0^E = 0$, there exists $T^{\ast}>0$ and a unique solution to the initial value problem
(\ref{EULER}) - (\ref{EULER-ic})
such that for all $0<T<T^{\ast}$ it holds,

$$\mathbf{u}^E\in C^k([0,T),W^{3-k,2}(\R^3;\R^3)), \,\Pi\in C^k([0,T),W^{3-k,2}(\R^3)),\, k=0,1,2,3$$

\begin{equation}\label{e2deuler}
\|\mathbf{u}^E\|_{W^{k,\infty}(0,T;W^{3-k,2}(\R^3;\R^3))} + \|\Pi\|_{W^{k,\infty}(0,T;W^{3-k,2}(\R^3))}\leq c(T,\|\mathbf{u}_0^E\|_{W^{3,2}(\R^3)}).
\end{equation}

\end{theorem}

\subsection{Main result} \label{MR}

Our main result rigorously characterises the aforementioned asymptotic regime of \eqref{mom} in the case of the following set of \emph{ill-prepared initial data}. Specifically,

\begin{equation}\label{eq:limitData-2}
    \varrho(0,\cdot)=\varrho_\eps^0=1+\eps \seps^0,\quad  \seps^0\in L^{\infty}(\R^3)\cap H^{1}(\R^3) \,\, \text{uniformly bounded}, \,\, \seps^0\rightarrow s^0 \quad \text{in} \,\, H^1(\R^3),
\end{equation}
    
 \begin{equation}
 \label{eq:limitData}
 \mathbf{u} (0,\cdot)=\mathbf{u}_\eps^0,\quad  \sqrt{\varrho_\eps^0}\mathbf{u}_\eps^0 \in L^{2}(\R^3), \quad   \sqrt{\varrho_\eps^0}(\mathbf{u}_\eps^0-\mathbf{u}^0)\rightarrow 0\quad \text{in} \quad L^2(\R^3).
\end{equation}

as $\varepsilon \to 0$.

\smallskip
Having recalled and collected all the preliminary notions we need in the sequel we are ready to state our main result.

\begin{theorem} \label{main}
    Suppose $\nu \to 0$ as $\varepsilon \to 0$. Assume the initial data $(\varrho_\varepsilon^0, \mathbf{u}_\varepsilon^0)$ to be of uniformly bounded
    finite energy and BD-entropy, i.e. there exists $C>0$ independent of $\varepsilon, \nu>0$ such that 
    \begin{equation}\label{eq: uniform bound data}
        E(\varrho_\eps^0,\mathbf{u}_\eps^0)\leq C, \quad   B(\varrho_\eps^0,\mathbf{u}_\eps^0)\leq C.
    \end{equation}
   and satisfy \eqref{eq:limitData-2}, \eqref{eq:limitData}. Let $(\varrho_\varepsilon, \mathbf{u}_\varepsilon, \mathcal{T}_\varepsilon)$ be a weak solution of (\ref{mom}) and $\mathbf{u}^E$ the unique solution to the initial value problem (\ref{EULER}) - (\ref{EULER-ic}) on $[0,T)\times \R^3$, $0<T<T^{\ast}$, with initial datum $\mathbf{u}_{0}^E = \mathbf{P}(\mathbf{u}^0) \in W^{3,2}(\mathbb{R}^3)$. Then, as $\varepsilon \to 0$,
\begin{equation}
\label{main-rho-conv}
\varrho_\varepsilon-1 \to 0 \ \ \text{in} \ \ L^\infty (0,T;L^2(\mathbb{R}^3)+L^\gamma(\mathbb{R}^3))\cap L^{\infty}(0,T;H^s(\R^3))
\end{equation}
for all $0<s<1$ and
\begin{equation}
\label{main-u-conv}
\sqrt{\varrho}_\varepsilon\mathbf{u}_\varepsilon \to \mathbf{u}^E
\ \ \text{in} \ \ L^2(0,T;L^2_{\mathrm{loc}}(\mathbb{R}^3),
\end{equation}
for any $0<T<T^{\ast}$.
\end{theorem} 

\begin{remark}
   Note that 
in \eqref{eq:limitData-2}, we do not require that the initial data for the density is well-prepared but allow for the presence of fast oscillating acoustic waves. The quantity $s^0$ will serve as initial data for the linear acoustic system in Section \ref{AWe} and subsequently in the relative entropy scheme in Section 5.
\end{remark}

\begin{remark}
We also provide the analogue of Theorem \ref{main} in the case $d=2$, the main difference being the weaker dispersive decay of the acoustic waves and the well-posedness properties of the $2D$-Euler equations \eqref{EULER}. Posed on $\R^2$, singularity formation in finite time is ruled out and we obtain global solutions. We refer to Section \ref{2D-c} for details.
\end{remark}

\section{Uniform bounds} \label{ID-UE}

Owing to \eqref{eq:EI} and \eqref{eq:BDI}, the weak solutions under consideration satisfy the following uniform bounds. To take the non-trivial farfield \eqref{farfield} into account,  we will repeatedly rely on the following fact: let $f\in L_{\loc}^1(\R^3)$ with $\nabla f\in L^2(\R^3)$, then there exists $c\in \R$ such that $f-c\in L^{6}(\R^3)$ (see \cite[Theorem 4.5.9]{Hormander} for the proof).
 
\begin{lemma}\label{lem:initial data}
Assume that the initial data $(\varrho_\varepsilon^0,\mathbf{u}_\eps^0)$ satisfies \eqref{eq: uniform bound data}. Then, the following hold:
\begin{enumerate}[(i)]
    \item $\varrho_\varepsilon^0-1\in H^1(\R^3)$ uniformly bounded and 
    \begin{equation*}
        \|\varrho_\varepsilon^0-1\|_{L^2(\R^3)}\leq C\begin{cases}
            \eps^{\frac{4}{6-\gamma}} \quad &1<\gamma<2\\
            \eps \quad &\gamma\geq 2,
        \end{cases}
    \end{equation*}
    \item $\sqrt{\varrho_\varepsilon^0}-1\in L^2(\R^3)$ uniformly bounded with
        \begin{equation*}
        \|\sqrt{\varrho_\varepsilon^0}-1\|_{L^2(\R^3)}\leq C\begin{cases}
            \eps^{\frac{4}{6-\gamma}} \quad &1<\gamma<2\\
            \eps \quad &\gamma\geq 2,
        \end{cases}
        \end{equation*}
        \item $\sqrt{\varrho_\varepsilon^0}\mathbf{u}_\eps^0\in L^2(\R^3)$ uniformly bounded.
    \item the initial momentum satisfies $\varrho_\varepsilon^0\mathbf{u}_\eps^0\in L^2(\R^3)+L^{\frac{3}{2}}(\R^3)$ uniformly bounded,
    \item the initial density fluctuations $\seps^0$ satisfy $\eps\seps^0\in H^1(\R^3)$ uniformly bounded and
    \begin{equation*}
        \|\seps^0\|_{L^2(\R^3)}\leq \begin{cases}
            C\eps^{\frac{\gamma-2}{6-\gamma}} \quad &1<\gamma<2,\\
            C \quad &\gamma \geq 2.
        \end{cases}
    \end{equation*}
\end{enumerate}
\end{lemma}

\begin{remark} \label{rem-E}
Note that the above uniform bounds are only based on the finite energy assumption and do not use additional information from the viscous dissipation and the BD-entropy inequality.
\end{remark}

\begin{proof}
The assumptions \eqref{eq: uniform bound data} yield that 
\begin{equation*}
    \nabla \varrho_\varepsilon^0\in L^2(\R^3), 
    \quad H(\varrho_\varepsilon^0)\in L^1(\R^3),
\end{equation*}
uniformly bounded and with $H$ defined as in \eqref{eq:energy} and thus non-negative. To prove (i), we note that by convexity of the function $s\rightarrow s^{\gamma}-1-\gamma(s-1)$ for $\gamma>1$ is possible to conclude that 
\begin{equation}\label{eq:Orlicz}
\int_{\R^3}\left|\varrho_\varepsilon^0-1\right|^{2} \mathbf{1}_{\{|\varrho_\varepsilon^0-1|\leq\frac{1}{2}\}}+\left|\varrho_\varepsilon^0-1\right|^{\gamma} \mathbf{1}_{|\varrho_\varepsilon^0-1|>\frac{1}{2}}\text{d}x\leq C\eps^2,
\end{equation}
see \cite{LM98} for details. If $\gamma\geq 2$, this yields
\begin{equation*}
\|\varrho_\varepsilon^0-1\|_{L^{2}(\R^3)}\leq C\eps
\end{equation*}
and in particular $\varrho_\varepsilon^0-1\in H^1(\R^3)$ uniformly bounded provided that $\gamma\geq 2$. Moreover, for any $\gamma>1$ 
it follows from $\nabla \varrho_\varepsilon^0\in L^2(\R^3)$ that there exists $c\in \R$ such that $\varrho_\varepsilon^0-c\in L^6(\R^3)$ uniformly bounded. The bound \eqref{eq:Orlicz} then entails that $\varrho_\varepsilon^0-1\in L^6(\R^3)$. For $1<\gamma<2$ it follows by interpolation that 
\begin{equation*}
    \|(\varrho_\varepsilon^0-1)\mathbf{1}_{|\varrho_\varepsilon^0-1|>\frac{1}{2}}\|_{L^2}\leq \|(\varrho_\varepsilon^0-1)\mathbf{1}_{|\varrho_\varepsilon^0-1|>\frac{1}{2}}\|_{L^\gamma}^{\theta
    }\|(\varrho_\varepsilon^0-1)\mathbf{1}_{|\varrho_\varepsilon^0-1|>\frac{1}{2}}\|_{L^6}^{1-\theta}\leq C\eps^{\frac{2}{\gamma
    }\theta}
\end{equation*}
with $\theta=\frac{2\gamma}{6-\gamma}$. 
Upon observing that 
$$\left|\sqrt{\varrho_\varepsilon^0}-1 \right| = \left|(1+\sqrt{\varrho_\varepsilon^0})^{-1}(\varrho_\varepsilon^0-1) \right| \leq \left| (\varrho_\varepsilon^0-1)\right|,$$
the desired bound (ii) is implied by (i).
While (iii) is a direct consequence of \eqref{eq: uniform bound data}, the statement (iv) follows from $\sqrt{\varrho_\varepsilon^0}-1\in H^1(\R^3)\hookrightarrow L^{6}(\R^3)$ and $\sqrt{\varrho_\varepsilon^0}\mathbf{u}_\eps^0\in L^2(\R^3)$ as
\[
\varrho_\varepsilon^0\mathbf{u}_\eps^0=\sqrt{\varrho_\varepsilon^0}\mathbf{u}_\eps^0+(\sqrt{\varrho_\varepsilon^0}-1)\sqrt{\varrho_\varepsilon^0}\mathbf{u}_\eps^0\in L^2(\R^3)+L^{\frac{3}{2}}(\R^3).
\]  
For the statement (v) it suffices to note that $\varrho_\varepsilon^0-1\in H^1(\R^3)$ uniformly bounded yields $\eps\seps^0\in L^2(\R^3)$ uniformly bounded. Moreover the $L^2$-decay rate for $\varrho_\varepsilon^0-1$ yields $\seps^0\in L^2(\R^3)$ uniformly bounded provided that $\gamma\geq 2$ and 
\begin{equation*}
\|\seps^0\|_{L^2(\R^3)}\leq C \eps^{\frac{\gamma-2}{6-\gamma}}
\end{equation*}
for $1<\gamma<2$.
\end{proof}
\noindent
Along the same lines of Lemma \ref{lem:initial data}, one infers the following uniform bounds for weak solutions to \eqref{mom} provided that \eqref{eq: uniform bound data} holds:
\begin{align} \label{bounds-1}
   &\|\varrho_\varepsilon-1 \|_{L_t^{\infty}H_x^1}\leq C, \quad  \|\varrho_\varepsilon-1\|_{L_t^{\infty}L_x^2}\leq C\begin{cases}
            \eps^{\frac{4}{6-\gamma}} \quad &1<\gamma<2,\\
            \eps \quad &\gamma\geq 2;
        \end{cases}\\   \label{bounds-2} &\|\sqrt{\varrho_\varepsilon}\mathbf{u}_\eps\|_{L_t^{\infty}L_x^2}\leq C, \quad   \|{\varrho_\varepsilon}\mathbf{u}_\varepsilon\|_{L_t^{\infty}L_x^2+L_x^{\frac{3}{2}}}\leq C.
\end{align}
By interpolation, one has that $\reps-1$ converges strongly to $0$ in $L^{\infty}(0,T;H^s(\R^3))$ for all $0\leq s<1$. In particular, for all $2\leq q< 6$ there exists $C>0$ such that 
\begin{equation}\label{eq: decay reps}
     \|\sqrt{\varrho_\varepsilon}-1\|_{L_t^{\infty}L_x^q}\leq \|\varrho_\varepsilon-1 \|_{L_t^{\infty}L_x^q}\leq C \eps^{\beta}
\end{equation}
where $\beta=\beta(q)=\frac{4}{5}\left(1-3\left(\frac{1}{2}-\frac{1}{q}\right)\right)$ by interpolation. 
\noindent
The viscous dissipation and the BD-entropy inequality further allow for the following bounds.

\begin{lemma} \label{lemm-BD}
    Assume that $(\varrho_\varepsilon,\ueps)$ is a finite energy weak solution to \eqref{mom} with initial data $(\varrho_\varepsilon^0,\ueps^0)$ satisfying \eqref{eq: uniform bound data}. Then there exits $C>0$ independent of $\eps$ such that
\begin{equation*}
\nu\|\sqrt{\varrho_\varepsilon}-1\|_{L_t^{\infty}H_x^1}\leq C, \quad \|\sqrt{\nu} \nabla \varrho_\varepsilon^{\frac{\gamma}{2}} \|_{L_t^{2}L_x^2}\leq C,
\end{equation*}
\begin{equation}\label{bounds-3}
\| \sqrt{\nu}\Delta \varrho_\varepsilon \|_{L_t^{2}L_x^2}\leq C, \quad
 \|\sqrt{\nu} \mathcal{S}(\mathbf{u})\|_{L_t^{2}L_x^2}\leq C, \quad \|\sqrt{\nu} \mathcal{A}(\mathbf{u})\|_{L_t^{2}L_x^2}\leq C.
\end{equation}
\end{lemma}

\begin{proof}
The bounds are immediate consequences of the energy and BD-entropy inequalities stated in Definition \ref{weak-def}.
\end{proof}

\section{Acoustic waves} \label{AWe}

The analysis of the incompressible limit requires a suitable control of the acoustic waves. Indeed, their presence is due to density fluctuations at finite Mach number and may create highly oscillating phenomena. Being the problem \eqref{mom} posed on $\R^d$, suitable dispersive estimates yield the decay of the acoustic waves in the low-Mach number limit. 
As will be detailed below, it is possible to obtain a fourth-order acoustic wave equation for the \emph{density fluctuations} 
$$\seps=\eps^{-1}(\varrho_\varepsilon-1)$$ 
by a suitable linearization of the system \eqref{mom} around the constant (incompressible) state, namely
\begin{equation}\label{eq:densityfluct}
    \partial_{tt}^2 \seps-\frac{1}{\eps^2}\Delta(1-2\eps^2\kappa^2\Delta) \seps=0.
\end{equation}
We notice that in the absence of the capillarity tensor $(\kappa=0)$, we have the usual dispersion relation for the classical wave equation. Here and due to the Korteweg tensor, this dispersion corresponds, at the leading order, to the one already observed in quantum fluids such as e.g. the QHD or quantum Navier-Stokes equations, see \cite{lars1, lars3}.
Denoting the sequence of momenta by $\mathbf{m}_\varepsilon=\varrho_\varepsilon\mathbf{u}_\varepsilon$ and linearizing  \eqref{mom} around the constant state $(\reps=1, \meps=0)$, one obtains the acoustic system
\begin{equation}\label{eq:acoustic1}
\begin{aligned}
        &\partial_t \seps+\frac{1}{\eps}\div \Q(\meps)=0,\\
        &\partial_t\Q(\mathbf{m}_\varepsilon)+\frac{\gamma}{\eps}\nabla\left( \seps-2\kappa^2\eps^2\Delta  \seps\right)=\Q( F_{\eps}),
\end{aligned}
\end{equation}
with
\begin{equation*}
    \Feps:=-\div\left(\varrho_\varepsilon\ueps\otimes \ueps\right)+2\nu\div\left(\varrho_\varepsilon\mathbb{D}\ueps\right) -\frac{(\gamma-1)}{\eps^2}\nabla H(\varrho_\varepsilon)+\eps^2G_{\eps}
\end{equation*}
and where $$\eps^{-2}\nabla p(\varrho_\varepsilon)=\frac{\gamma}{\eps}\nabla \seps+\frac{(\gamma-1)}{\eps^2}\nabla H(\varrho_\varepsilon).$$
In particular, 
for the dispersive stress tensor $\div \mathbb{K}$, as defined in \eqref{K}, we obtain the identity
\begin{align*}
    \div\mathbb{K}&=\div\left((1+\eps \seps)\div(\eps\nabla \seps)\mathbb{I}+\frac{1}{2}\left|\eps\nabla \seps|^2\right)\mathbb{I}\right)-\div\left(\eps\nabla \seps\otimes \eps\nabla \seps\right)\\&=\eps\nabla\Delta \seps+\eps^2 G_{\eps}\\
\end{align*}
with 
$$G_{\eps}:=\seps \nabla\Delta \seps.$$

While we discarded all terms in $\div\mathbb K$ being multiplied by $\eps^2$ into $\Feps$, we keep the leading order linear term $2\eps\kappa^2\nabla\Delta \seps$ which, as already mentioned, alters the classical dispersion relation. Consequently, the decay of $\sigma_\varepsilon$ in space-time norms  can be shown to be stronger compared to the classical wave equations, see Proposition \ref{prop:Strichartz} below and \cite[Remark 4.10]{lars2}. The method of exploiting altered dispersion relations through accurate Strichartz estimates has also been employed e.g. in the analysis of the quasi-neutral limit for Navier-Stokes-Korteweg system \cite{donatella-piero}.

We sketch the derivation of the respective Strichartz estimates and refer \cite[Section 4]{lars1} to for further details. We symmetrize \eqref{eq:acoustic1} by applying 
\begin{equation}\label{eq:sym-op}
    \widetilde{\seps}=(1-2\eps^2\kappa^2\Delta)^{\frac{1}{2}} \seps, \qquad \tmeps=(-\Delta)^{-\frac{1}{2}}\div\meps
\end{equation}
to obtain 
\begin{equation}\label{eq:diagonal}
\begin{aligned}
&\partial_t \widetilde{\seps}+\frac{1}{\eps}(-\Delta)^{\frac{1}{2}}(1-2\eps^2\kappa^2\Delta)^{\frac{1}{2}}\tmeps=0,\\
&\partial_t\tmeps-\frac{1}{\eps}(-\Delta)^{\frac{1}{2}}(1-2\eps^2\kappa^2\Delta)^{\frac{1}{2}} \widetilde{\seps}=\tFeps,
\end{aligned}
\end{equation}
where 
\begin{equation*}
    \tFeps=(-\Delta)^{-\frac{1}{2}}\div\Q\left(\Feps\right).
\end{equation*}
As described in \cite{lars1}, the system \eqref{eq:diagonal} can be characterized by means of the linear semi-group operator $\eith$ where $H_{\eps}$ is defined via the Fourier multiplier 
\begin{equation} \label{F-mult}
    \phie(|\xi|)=\frac{|\xi|}{\eps}\sqrt{1+2\eps^2\kappa^2|\xi|^2}. 
\end{equation}
The dispersion given by \eqref{F-mult} mimics wave-like behavior for low frequencies $|\xi|<\eps^{-1}$ while it is Schr\"odinger like for high-frequencies $|\xi|>\eps^{-1}$. A stationary phase argument then leads to the desired Strichartz (dispersive) estimates \cite[Section 4]{lars1}. 

In order to recall the dispersive estimates suitable for our analysis, we start saying that a pair of Lebesgue exponents $(p, q)$ is called \emph{``Schr\"odinger admissible"} if 
$$2\leq p, q\leq\infty \quad  \text{and} \quad \frac2p+\frac3q=\frac 32.$$ 
The following proposition holds.
\begin{proposition}[{\cite[Proposition 4.2]{lars1}}]\label{prop:Strichartz}
Fix $\eps>0$ and $s\in \R$. There exists a constant $C>0$ independent from $T$ and $\eps$ such that for any $(p,q)$ admissible pair and any $\alpha_0\in [0,\frac{1}{2}(\frac12-\frac1q)]$, the following hold true,
\begin{equation}\label{eq:Strichartz-f1}
\|\eith f\|_{L^p(0,T;W^{s-\alpha_0,q}(\R^3))}\leq C \eps^{\alpha_0}\|f\|_{H^{s}(\R^3)}.
\end{equation}
\end{proposition}

For the purpose of the inviscid incompressible limit (see Section \ref{IIL}) it turns out to be sufficient to consider the acoustic system system \eqref{eq:acoustic1} in its homogeneous form
\begin{equation}\label{eq:AW-linear}
\begin{aligned}
        &\partial_t s_{\eps}+\frac{1}{\eps}\div\nabla \Phi_{\eps}=0,\\
        &\partial_t\nabla\Phi_{\eps}+\frac{\gamma}{\eps}\nabla\left( s_{\eps}-2\kappa^2\eps^2\Delta  s_{\eps}\right)=0,
\end{aligned}
\end{equation}
with initial data $(s^0, \nabla\Phi_{0})$, where we denoted by $\Phi_{\eps}$ the acoustic potential, i.e. $\nabla \Phi_{\eps} =\Q(\meps) $. 
For the sake of clarity, we choose a distinct notation for the solution $s_\eps$ of the linear homogeneous equation \eqref{eq:AW-linear} and the solution $\seps$ to \eqref{eq:acoustic1}. 
In particular, we will see in Section \ref{IIL} that the initial data $(s^0, \nabla\Phi_{0})$ considered will be given by the data $(s^0, \mathbf{Q}(\mathbf{u}_0))$ in \eqref{eq:limitData}, \eqref{eq:limitData-2}. 
properly regularized.

\begin{remark}
    Note that the fourth-order equation \eqref{eq:densityfluct} formally follows from \eqref{eq:AW-linear} by applying the differential operators $\partial_t$ and $\div$ to the first and second equation respectively.
\end{remark}

From Proposition \ref{prop:Strichartz} and combing the transformation \eqref{eq:sym-op} with the Strichartz estimates of Proposition \ref{prop:Strichartz}, one derives the following bounds for solutions to \eqref{eq:AW-linear}. we refer the reader to \cite[Proposition 4.5]{lars1} for details.
\begin{proposition} \label{prop-AW-dec}
    Let $s\in \R$, $(s^0,\mathbf{m}^0)$ be such that $s^0\in H^s(\R^3)$, $\eps\nabla s^0\in H^s(\R^3)$ and $\Q (\mathbf{m}^0)\in H^s(\R^3)$ and denote by $(s_{\eps},\meps)$ the unique solution to the homogeneous system \eqref{eq:AW-linear} with initial data $(s^0, \Q(\mathbf{m}^0))$. Then, for all $T>0$, all Schr\"odinger admissible pair $(p,q)$ and any $\alpha\in [0,\frac{1}{2}(\frac12-\frac1q)]$ the following holds true
    \begin{equation*}
        \|(s_\eps,\Q(\meps))\|_{L^{p}(0,T;W^{s,q})\R^3)} 
    \end{equation*}
    \begin{equation*}
        \leq C\eps^{\alpha}\left(\|s^0\|_{H^{s+\alpha}(\R^3)}+\|\eps\nabla s^0\|_{H^{s+\alpha}(\R^3)}+\|\Q(\mathbf{m}^0)\|_{H^{s+\alpha}(\R^3)}\right).
        \end{equation*}
\end{proposition}

Note that $(s_\eps, \Q(\meps))$ can be bounded in terms of $( \widetilde{\seps}, \tmeps)$ in $L^p(0,T;W^{s,q}(\R^3))$ in view of \eqref{eq:sym-op}, while one has 
\begin{equation*}
    \left\|(\tilde{s}_\eps^0, \tilde{m}^0) \right\|_{H^{s}(\R^3)} \leq \|s^0\|_{H^{s+\alpha}(\R^3)}+\|\eps\nabla s^0\|_{H^{s+\alpha}(\R^3)}+\|\Q(\mathbf{m}^0)\|_{H^{s}(\R^3)}
\end{equation*}
which leads to the above estimate. 

\section{Relative energy inequality} \label{REI}

In the same spirit of \cite{BCD}, we derive the relative energy inequality related to the Navier-Stokes-Korteweg system (\ref{mom}).
For $\tau\in[0,T]$, we introduce the following relative energy functional
$$
\mathcal{E}(\varrho, \mathbf{u}\,|\, r, \mathbf{U}) (\tau) =  \int_{\R^3} \frac{1}{2} |\sqrt{\varrho}\mathbf{u}
- \sqrt{\varrho}\,\mathbf{U}|^2  (\tau,\cdot)  dx
+
\kappa^2 \int_{\R^3} |\nabla \varrho
- \nabla r
|^2  (\tau,\cdot)  dx
$$
$$
+ \frac{1}{\varepsilon^2} \int_{\R^3} \left[H(\varrho) - H(r)
- H'(r)(\varrho - r)\right](\tau,\cdot)  dx,
$$
where  $(\varrho,\mathbf{u})$ is a weak solution to $\eqref{mom}$ and
$(r,\mathbf{U})$ is a pair of smooth (arbitrary) test functions. 
In the following,
we simply write $\mathcal{E}(\tau)$ in place of $\mathcal{E}(\varrho, \mathbf{u}\,|\, r, \mathbf{U}) (\tau)$.
First, thanks to the energy inequality (\ref{eq:EI}), we infer that 
\begin{equation} \label{step-1}
\begin{aligned}
\mathcal{E}(s)\big|_{s=0}^{s=\tau} 
&\leq  \int_{\R^3} \left(\frac{1}{2} \varrho |\mathbf{U}|^2 - \varrho \mathbf{u} 
\cdot \mathbf{U} \right)(s,\cdot) dx\big|_{s=0}^{s=\tau}
+ \kappa^2
\int_{\R^3} \left(|\nabla r|^2 - 2 \nabla \varrho \cdot \nabla r \right)(s,\cdot) dx\big|_{s=0}^{s=\tau}
\\
&- \frac{1}{\varepsilon^2} \int_{\R^3} \left( H(r) + H'(r)(\varrho-r) \right) (s,\cdot)dx \big|_{s=0}^{s=\tau}\\
&\quad- \nu\int_0^\tau\int_{\R^3}\left|\mathcal{S}(\mathbf{u})\right|^2dxdt.
\end{aligned}
\end{equation}
Next, we test the continuity equation by $\frac{1}{2}|\mathbf{U}|^2$ and $2\kappa^2 \Delta r$, and the momentum equation by $\mathbf{U}$,
to get
$$
\int_0^\tau \int_{\R^3}\frac{1}{2}\varrho \partial_t |\mathbf{U}|^2\, dx dt 
+\int_0^\tau \int_{\R^3}\frac{1}{2} \varrho \mathbf{u} \cdot \nabla |\mathbf{U}|^2 \,dx dt  
= \int_{\R^3} \frac{1}{2} \varrho |\mathbf{U}|^2 (s,\cdot) dx\big|_{s=0}^{s=\tau},
$$
$$
-2\kappa^2\int_0^\tau \int_{\R^3}\varrho \partial_t \Delta r\, dx dt 
-2\kappa^2\int_0^\tau \int_{\R^3} \varrho \mathbf{u} \cdot \nabla \Delta r \,dx dt  
= 2\kappa^2\int_{\R^3} 
\nabla \varrho \cdot \nabla r(s,\cdot) dx\big|_{s=0}^{s=\tau}
$$
and
\begin{gather*}
-\int_0^\tau \int_{\R^3} \varrho \mathbf{u}\cdot \partial_t \mathbf{U} \, dxdt
-\int_0^\tau \int_{\R^3} (\sqrt{\varrho} \mathbf{u} 
\otimes \sqrt{\varrho} \mathbf{u}) : \nabla \mathbf{U} \, dxdt 
-\frac{1}{\varepsilon^2}\int_0^\tau \int_{\R^3} p(\varrho) \textrm{div}\,\mathbf{U} \,dxdt \\\noalign{\vskip0.5mm}
-2\kappa^2 \int_0^{\tau} \int_{\R^3} 
\left( 
\nabla\varrho \cdot \nabla( \varrho \textrm{div}\mathbf{U}) 
- \frac{1}{2} |\nabla \varrho|^2 \textrm{div}\mathbf{U}
+ \nabla \varrho \otimes \nabla \varrho : \nabla \mathbf{U}
\right)dxdt
\\\noalign{\vskip0.5mm}
+2\nu\int_0^\tau \int_{\R^3} \sqrt{\varrho}
\mathcal{S}(\mathbf{u})
: \nabla \mathbf{U} dxdt 
= -\int_{\R^3} (\varrho \mathbf{u} \cdot \mathbf{U})(s,\cdot) dx\big|_{s=0}^{s=\tau}.
\end{gather*}
Moreover, we have
$$
\kappa^2\int_{\R^3} 
|\nabla r(s,\cdot)|^2 dx\big|_{s=0}^{s=\tau}
=
-2\kappa^2 \int_{0}^\tau\int_{\R^3}
\partial_t r
\Delta r dxdt.
$$
Then, using the continuity equation, we have
\begin{equation}  \label{test-H}
 \hspace{-0.3 cm} \begin{aligned}
  & \int_0^T\!\! \int_{\R^3} \left[ \partial_t \left( H(r)
 + H'(r)(\varrho - r) \right) \right] dxdt
\\
& \,\,= \int_0^T \!\!\int_{\R^3} \left[ H'(r)\partial_t r + \partial_t (H'(r))(\varrho - r)
  + H'(r)\partial_t \varrho - H'(r) \partial_t r
  \right] dxdt
\\
 &\,\, = \int_0^T\!\! \int_{\R^3} \left[ \partial_t (H'(r))(\varrho - r)
  - H'(r)\div_x(\varrho \mathbf{u})  \right] dxdt
\\
& \,\, = \int_0^T\!\! \int_{\R^3}  \left[  \partial_t  (H'(r))(\varrho - r) + \varrho
    \mathbf{u} \cdot \nabla_x (H'(r))  \right]  dxdt.
\end{aligned}
\end{equation}
Therefore, from (\ref{step-1}), we obtain
\begin{equation}
\label{step-2}
\begin{aligned}
\mathcal{E}(s)\big|_{s=0}^{s=\tau} 
&\leq  \int_0^\tau \int_{\R^3}  \frac{1}{2}\varrho \partial_t |\mathbf{U}|^2\, dx dt + \int_0^\tau \int_{\R^3}
\frac{1}{2}\varrho \mathbf{u} \cdot \nabla |\mathbf{U}|^2 \,dx dt
\\&\quad - \int_0^\tau \int_{\R^3} \varrho \mathbf{u}\cdot \partial_t \mathbf{U} \, dxdt
-\int_0^\tau \int_{\R^3} (\sqrt{\varrho} \mathbf{u} \otimes\sqrt{\varrho} \mathbf{u}) : \nabla \mathbf{U} \, dxdt 
\\&\quad
+2\kappa^2 \int_0^{\tau} \int_{\R^3}
(\varrho\partial_t \Delta r + \varrho\mathbf{u} \cdot \nabla\Delta r
-\partial_t r \Delta r
)dxdt
\\&\quad 
-2\kappa^2 \int_0^{\tau} \int_{\R^3} 
\left( 
\nabla\varrho \cdot \nabla( \varrho \textrm{div}\mathbf{U}) 
- \frac{1}{2} |\nabla \varrho|^2 \textrm{div}\mathbf{U}
+ \nabla \varrho \otimes \nabla \varrho : \nabla \mathbf{U}
\right)dxdt
\\&\quad 
+2\nu\int_0^\tau \int_{\R^3} \sqrt{\varrho}
\mathcal{S}(\mathbf{u})
: \nabla \mathbf{U} dxdt 
\\&\quad 
- \frac{1}{\varepsilon^2} \int_0^\tau\!\! \int_{\R^3}  \left[  \partial_t  (H'(r))(\varrho - r) 
+ \varrho\mathbf{u} \cdot \nabla (H'(r)) +p(\varrho) \mbox{div} \mathbf{U} \right]  dxdt
\\&\quad 
-\nu\int_0^\tau\int_{\R^3}\left|\mathcal{S}(\mathbf{u})\right|^2dxdt. 
\end{aligned}
\end{equation}
Rearranging the relative energy inequality above, we obtain
\begin{equation}
\label{step-3}
\begin{aligned}
\mathcal{E}(s)\big|_{s=0}^{s=\tau} 
+&\nu\int_0^\tau\int_{\R^3}\left|\mathcal{S}(\mathbf{u})\right|^2dxdt
\\& \leq  \int_0^\tau \int_{\R^3}
\left[\partial_t \mathbf{U} \cdot (\varrho\mathbf{U} - \varrho\mathbf{u}) 
+ (\sqrt{\varrho}\mathbf{u} \cdot \nabla) \mathbf{U} \cdot (\sqrt{\varrho}\mathbf{U} - \sqrt{\varrho}\mathbf{u})\right]\,dx dt 
\\&
+2\kappa^2 \int_0^{\tau} \int_{\R^3}
(\varrho\partial_t \Delta r + \varrho\mathbf{u} \cdot \nabla\Delta r
-\partial_t r \Delta r
)dxdt
\\& 
-2\kappa^2 \int_0^{\tau} \int_{\R^3} 
\left( 
\nabla\varrho \cdot \nabla( \varrho \textrm{div}\mathbf{U}) 
- \frac{1}{2} |\nabla \varrho|^2 \textrm{div}\mathbf{U}
+ \nabla \varrho \otimes \nabla \varrho : \nabla \mathbf{U}
\right)dxdt
\\& 
+2\nu\int_0^\tau \int_{\R^3} \sqrt{\varrho}
\mathcal{S}(\mathbf{u})
: \nabla \mathbf{U} dxdt 
\\& 
- \frac{1}{\varepsilon^2} \int_0^\tau\!\! \int_{\R^3}  \left[  \partial_t  (H'(r))(\varrho - r) 
+ \varrho\mathbf{u} \cdot \nabla (H'(r)) +p(\varrho) \mbox{div} \mathbf{U} \right]  dxdt
\\& = I_1 + ... + I_5.
\end{aligned}
\end{equation}
Apart from the presence of the Korteweg terms, relation \eqref{step-3} is reminiscent of the one derived in \cite{BCD} (see relation (4.7)).

\section{Inviscid incompressible limit} \label{IIL}

For a fixed $\delta > 0$, we consider the following ansatz
\begin{equation} \label{ansatz}
    r=1+\varepsilon s_{\eps,\delta}, \ \ \ \ \mathbf{U} = \mathbf{u}^E + \nabla \Phi_{\eps,\delta}
\end{equation}
in the relative energy inequality (\ref{step-3}). Precisely, $\mathbf{u}^{E}$ is the solution of the target system  \eqref{EULER}  with initial datum $\mathbf{u}^E_0=\mathbf{P}(\mathbf{u}_{0})$ and $s_{\eps,\delta}$, 
$\nabla\Phi_{\eps,\delta}$ are the solution of the acoustic system \eqref{eq:AW-linear} with the regularized initial data 
\begin{equation}\label{eq:data}
    s_{\delta}^0=s^{0}\ast \eta_\delta, \qquad \nabla\Phi_{\delta}^0=\nabla\Phi^0\ast\eta_\delta, \qquad \Q (\mathbf{u}^0)=\nabla \Phi^0,
\end{equation}
for $\eta_\delta\in C_c^{\infty}(\R^3)$ a standard  mollifier and $(s^0, \mathbf{u}^0)$ as in \eqref{eq:limitData-2} and  \eqref{eq:limitData}.
The above choice is motivated by the fact that we consider the (arbitrary) smooth pair $(r,\mathbf{U})$ as the sum of the incompressible Euler equations and the contribution from acoustic waves. The latter disappear in the low-Mach number limit as the next Corollary shows. By means of the dispersive estimates  of Proposition \ref{prop-AW-dec}, we infer the following decay in 
$\eps$ in Strichartz norms.
\begin{corollary} \label{cor-AW-dec}
Let $s\in \R$, and $(s_{\delta}^0, \nabla \Phi_{\delta}^0)$ be as in \eqref{eq:data}. For any Schr\"odinger-admissible pair $(p,q)$ with $q>2$ and $\alpha,\delta>0$ sufficiently small and $T>0$ there exists $C>0$ only depending on the constant in \eqref{eq: uniform bound data} and $\delta>0$ such that the unique solution $(s_{\eps,\delta}, \nabla \Phi_{\eps,\delta})$ to the \eqref{eq:AW-linear} satisfies
\begin{equation}\label{eq:Strichartz final}
        \|(s_{\eps,\delta}, \nabla \Phi_{\eps,\delta})\|_{L^{p}(0,T;W^{s,q})\R^3)}\leq C\eps^{\alpha}.
\end{equation}
\end{corollary}

Note, that the solution $(s_{\eps,\delta},\nabla\Phi_{\eps,\delta})$ satisfies the energy conservation
\begin{equation} \label{en-cons}
    \|(s_{\eps,\delta}, \nabla \Phi_{\eps,\delta})\|_{L^{\infty}(\R,L^2(\R^3))}=\|( s_{\delta}^0,\nabla\Phi_{\delta}^0)\|_{L^2(\R^3)},
\end{equation}
and similarly for all $s\in \R$ it holds
\begin{equation} \label{Hs-cons}
    \|(s_{\eps,\delta}, \nabla \Phi_{\eps,\delta})\|_{L^{\infty}(\R,H^s(\R^3))}=\|( s_{\delta}^0,\nabla\Phi_{\delta}^0)\|_{H^s(\R^3)}
\end{equation}
due to the fact that $\eith$ constitutes a unitary semigroup operator.

In the following, we will consider the ansatz \eqref{ansatz} and we will handle the terms $I_1,...,I_5$ in (\ref{step-3}) in order to bound the relative entropy functional. We start with the relative energy functional for the initial data, whose required convergence (see \eqref{eq:limitData}, \eqref{eq:limitData-2}) is crucial in order to prove Theorem \ref{main} in terms of Gronwall's type arguments (see Section \ref{proof-closing} below).

\subsection{Initial data} \label{id}
From  the assumptions on the initial data $(\varrho_\eps^0,\mathbf{u}_\eps^0)$ in \eqref{eq:limitData-2}-\eqref{eq:limitData} and $(r,\mathbf{U})$ in \eqref{eq:data} respectively, we infer that
	$$
	\mathcal{E}(\varrho_\eps, \mathbf{u}_\eps\,|\, r, \mathbf{U}) (0)
	$$
	$$
	= \int_{\mathbb{R}^3}\frac{1}{2}\left|\sqrt{\varrho_\varepsilon^0}(\mathbf{u}_\varepsilon^0-\mathbf{u}_{0})\right|^{2} dx
	+\kappa^2\varepsilon^2\int_{\mathbb{R}^3}\left|\nabla \seps^0 - \nabla s_{\delta}^0 \right|^{2} dx
	$$
	\begin{equation} \label{initial data conv}
	+\int_{\mathbb{R}^3}\frac{1}{\varepsilon^{2}}\left[H\left(1+\eps\seps^0\right)-\varepsilon H^{\prime}\left(1+\eps s_{\delta}^0\right)\left(\seps^0-s_{\delta}^0\right)-H\left(1+\eps s_{\delta}^0\right)\right]dx.
	\end{equation}
	Now, setting $a=1+\eps\seps^0$ and $b=1+\eps s_{\varepsilon,\delta}$, we observe that there exists $\eps_0>0$ and $C>0$ independent of $\eps$ such that for all $\eps\in(0,\eps_0)$ it holds
	$$
	H(a)=H(b)+H^{\prime}(b)(a-b)+\frac{1}{2}H^{\prime\prime}(\xi)(a-b)^{2},\;\;\xi\in\left(a,b\right),
	$$
	$$
	\left|H(a)-H^{\prime}(b)(a-b)-H(b)\right|\leq C\left|a-b\right|^{2}.
	$$
    Consequently, we have
	$$
	\int_{\mathbb{R}^3}\frac{1}{\varepsilon^{2}}\left[H\left(1+\eps\seps^0\right)-\varepsilon H^{\prime}\left(1+\eps s_{\delta}^0\right)\left(\seps^0-s_{\delta}^0\right)-H\left(1+\eps s_{\delta}^0\right)\right]dx
	$$
	\begin{equation} \label{initial data conv2}
	\leq C\int_{\mathbb{R}^3}\frac{1}{\varepsilon^{2}}\left|\varepsilon\left(\seps^0-s_{\delta}^0\right)\right|^{2}dx= C\left\Vert\seps^0-s_{\delta}^0\right\Vert _{L^{2}(\mathbb{R}^{3})}^2.
	\end{equation}
	It follows
	$$
	\mathcal{E}(\varrho_\eps, \mathbf{u}_\eps\,|\, r, \mathbf{U}) 
	$$
	\begin{equation} \label{conv-id}
	\leq C\left\Vert\sqrt{\varrho_\varepsilon^0}(\mathbf{u}_\varepsilon^0-\mathbf{u}_{0})\right\Vert _{L^{2}(\mathbb{R}^{3};\mathbb{R}^{3})}^2
    + C(\kappa^2)\varepsilon^2 
    \left\Vert \nabla \seps^0-\nabla s_{\delta}^0\right\Vert _{L^{2}(\mathbb{R}^{3}; \mathbb{R}^{3})}^2+C\left\Vert\seps^0-s_{\delta}^0\right\Vert _{L^{2}(\mathbb{R}^{3})}^2.
	\end{equation}

\subsection{Convective terms} \label{ct}
Next, we provide bounds for the convective terms corresponding to $I_1$ in \eqref{step-3}. In the following we drop the subscript $\eps$ in $(\reps,\ueps)$ for the sake of a concise notation. Specifically, we have
$$
I_1
=
\int_0^\tau \int_{\R^3}
\left[\partial_t \mathbf{U} \cdot (\varrho\mathbf{U} - \varrho\mathbf{u}) 
+ (\sqrt{\varrho}\mathbf{u} \cdot \nabla) \mathbf{U} \cdot (\sqrt{\varrho}\mathbf{U} - \sqrt{\varrho}\mathbf{u})\right]\,dx dt 
$$

$$
= \int_0^\tau \int_{\R^3}
\left(\partial_t \mathbf{U} + \mathbf{U} \cdot \nabla \mathbf{U} \right)\cdot (\varrho\mathbf{U} - \varrho\mathbf{u}) 
dx dt 
$$

$$
+ \int_0^\tau \int_{\R^3}
\nabla \mathbf{U} \cdot (\sqrt{\varrho}\mathbf{U} - \sqrt{\varrho}\mathbf{u}) (\sqrt{\varrho}\mathbf{U} - \sqrt{\varrho}\mathbf{u})\,dx dt = I_1^{(1)} + I_1^{(2)}.
$$
As $\nabla U\in L_{t,x}^\infty$, the term $I_1^{(2)}$ is controlled by 
$$
|I_1^{(2)}|
\leq
C\int_0^\tau \mathcal{E}(\cdot,t)dt 
$$ 
while $I_1^{(1)}$ can be written as follows

$$
I_1^{(1)} = \int_0^\tau \int_{\R^3} 
(\varrho\mathbf{U} - \varrho\mathbf{u}) \cdot
(\partial_t \mathbf{u}^E + \mathbf{u}^E \cdot \nabla \mathbf{u}^E)\; dxdt
$$

$$
+ \int_0^\tau \int_{\R^3} 
(\varrho\mathbf{U} - \varrho\mathbf{u}) \cdot \partial_t \nabla \Phi_{\eps,\delta}\; dxdt
+ \int_0^\tau \int_{\R^3} 
(\varrho\mathbf{U} - \varrho\mathbf{u}) \otimes \nabla \Phi_{\eps,\delta} : \nabla \mathbf{u}^E \; dxdt
$$

\begin{equation} \label{I_1-split}
+ \int_0^\tau \int_{\R^3} 
(\varrho\mathbf{U} - \varrho\mathbf{u}) \otimes \mathbf{u}^E : \nabla^2 \Phi_{\eps,\delta}\; dxdt
+ \frac{1}{2}\int_0^\tau \int_{\R^3} 
(\varrho\mathbf{U} - \varrho\mathbf{u}) \cdot
\nabla |\nabla \Phi_{\eps,\delta}|^2 \; dxdt.
\end{equation}
Now, we analyze the first term in (\ref{I_1-split}). We have
$$
\int_0^\tau \int_{\R^3} 
(\varrho\mathbf{U} - \varrho\mathbf{u}) \cdot
(\partial_t \mathbf{u}^E + \mathbf{u}^E \cdot \nabla \mathbf{u}^E)\; dxdt = \int_0^\tau \int_{\R^3} \varrho\mathbf{u}\cdot \nabla \Pi \; dxdt
- \int_0^\tau \int_{\R^3} \varrho\mathbf{U}\cdot \nabla \Pi \; dxdt.
$$
In view of \eqref{bounds-2}, there exists $\overline{\mathbf{u}}\in L^{\infty}(0,T;L^2(\R^3))$ such that $\sqrt{\varrho}\mathbf{u} \rightharpoonup^{\ast} \overline{\mathbf{u}}$ in $ L_t^{\infty}L_x^2$. Together with \eqref{eq: decay reps}, it follows that 
$$
{\varrho}\mathbf{u}=(\sqrt{\rho}-1)\sqrt{\rho}u+\sqrt{\rho}u \to \overline{\mathbf{u}}
\ \ \ \ 
\text{weakly-(*) in} \ L_t^{\infty}L_x^2+L_x^{\frac{3}{2}-}.
$$
In particular, one has 
we deduce that 
$$
\int_0^{\tau} 
\int_{\R^3}
\varrho \mathbf{u} \cdot \nabla \Pi
\;
dx dt 
\to
\int_0^{\tau} 
\int_{\R^3}
\overline{\mathbf{u}} \cdot \nabla \Pi 
\;
dx dt.
$$
We consider the weak formulation of the continuity equation
$$
\int_0^{\tau} 
\int_{\R^3}
\left(
\varrho \partial_t \varphi + \varrho \mathbf{u} \cdot \nabla \varphi 
\right)
dx dt = \int_{\R^3}
\varrho (\tau,\cdot) \varphi (\tau,\cdot) dx - \int_{\R^3}
\varrho (0,\cdot) \varphi (0,\cdot) dx
$$
and for $\Pi$ as test-function to obtain
$$
\int_0^{\tau} 
\int_{\R^3}
\varrho{\mathbf{u}} \cdot \nabla \Pi
\;
dx dt \rightarrow\int_0^{\tau} 
\int_{\R^3}
\overline{\mathbf{u}} \cdot \nabla \Pi
\;
dx dt
$$
and
$$
\left[
\int_{\R^3} \varrho\Pi\; dx |_0^\tau - \int_0^{\tau} 
\int_{\R^3} \varrho \partial_t \Pi\;dxdt
\right]
\to 0 \ \ \ \ \text{as} \ \ \ \ \varepsilon\to0
$$
thanks to Lemma \ref{lem:initial data} and (\ref{bounds-1}). 
Consequently,
$$
\int_0^{\tau} 
\int_{\R^3}
\varrho \mathbf{u} \cdot \nabla \Pi
\;
dx dt
\to 0 \ \ \ \ \text{as} \ \ \ \ \varepsilon\to0.
$$
Next,
$$
\int_0^\tau \int_{\R^3} \varrho\mathbf{U}\cdot \nabla \Pi \; dxdt = 
\int_0^\tau \int_{\R^3} (\varrho-1)\mathbf{U}\cdot \nabla \Pi \; dxdt + 
\int_0^\tau \int_{\R^3} \mathbf{U}\cdot \nabla \Pi \; dxdt.
$$
With similar arguments as above
$$
\int_0^\tau \int_{\R^3} (\varrho-1)\mathbf{U}\cdot \nabla \Pi \; dxdt \to 0 \ \ \ \ \text{as} \ \ \ \ \varepsilon\to0,
$$
while
$$
\int_0^\tau \int_{\R^3} \mathbf{U}\cdot \nabla \Pi \; dxdt
= 
\int_0^\tau \int_{\R^3} \mathbf{u}^E\cdot \nabla \Pi \; dxdt
+
\int_0^\tau \int_{\R^3} \nabla \Phi_{\eps,\delta}\cdot \nabla \Pi \; dxdt.
$$
Performing integration by parts,
$$
\int_0^\tau \int_{\R^3} \text{div}\mathbf{u}^E\cdot \Pi \; dxdt = 0
$$
thanks to the incompressibility condition $\text{div}_{x}\mathbf{u}^E =0$. For the other term, using again integration by parts and the acoustic equation \eqref{eq:AW-linear}, we have
$$
\int_0^\tau \int_{\R^3} \nabla \Phi_{\eps,\delta}\cdot \nabla \Pi \; dxdt
=
- \int_0^\tau \int_{\R^3} \Delta \Phi_{\eps,\delta} \Pi \; dxdt= \varepsilon \int_0^\tau \int_{\R^3}
\partial_t s_{\varepsilon,\delta} \; \Pi \;
dxdt
$$

$$
= 
\left[
\varepsilon \int_{\R^3} s_{\varepsilon,\delta} \; \Pi \; dx |_0^\tau
- \varepsilon \int_0^\tau \int_{\R^3} 
s_{\varepsilon,\delta} \; \partial_t \Pi \; dxdt
\right]\to 0 \ \ \ \ \text{as} \ \ \ \ \varepsilon\to0.
$$
Now, we analyze the second term in (\ref{I_1-split}). We have
$$
\int_0^\tau \int_{\R^3} 
(\varrho\mathbf{U} - \varrho\mathbf{u}) \cdot \partial_t \nabla \Phi_{\eps,\delta}\; dxdt
= 
\int_0^\tau \int_{\R^3} 
\varrho\mathbf{u}^E \cdot \partial_t \nabla \Phi_{\eps,\delta}\; dxdt
$$
$$
- \int_0^\tau \int_{\R^3} 
\varrho\mathbf{u} \cdot \partial_t \nabla \Phi_{\eps,\delta}\; dxdt
+ \frac{1}{2} \int_0^\tau \int_{\R^3} \varrho\partial_t|\nabla \Phi_{\eps,\delta}|^2 \; dxdt
$$
$$
= \int_0^\tau \int_{\R^3} 
\varrho\mathbf{u}^E \cdot \partial_t \nabla \Phi_{\eps,\delta}\; dxdt
+ I^{I} + I^{II},
$$
where
\begin{equation}
\label{I-II}
I^{I} = - \int_0^\tau \int_{\R^3} 
\varrho\mathbf{u} \cdot \partial_t \nabla \Phi_{\eps,\delta}\; dxdt, \ \ \ \ I^{II} = \frac{1}{2} \int_0^\tau \int_{\R^3} \varrho\partial_t|\nabla \Phi_{\eps,\delta}|^2 \; dxdt,
\end{equation}
The terms $I^{I}$ and $I^{II}$ will be handled later in combination with other terms.
For the remaining one, we have
$$
\int_0^\tau \int_{\R^3} 
\varrho\mathbf{u}^E \cdot \partial_t \nabla \Phi_{\eps,\delta}\; dxdt
=
\int_0^\tau \int_{\R^3} 
(\varrho - 1)\mathbf{u}^E \cdot \partial_t \nabla \Phi_{\eps,\delta}\; dxdt
+
\int_0^\tau \int_{\R^3} 
\mathbf{u}^E \cdot \partial_t \nabla \Phi_{\eps,\delta}\; dxdt
$$
where, thanks to $\text{div}_{x}\mathbf{u}^E =0$, the second term
$$
\int_0^\tau \int_{\R^3} 
\text{div}_{x}\mathbf{u}^E \cdot \partial_t \Phi_{\eps,\delta}\; dxdt = 0.
$$
Using the acoustic equation \eqref{eq:AW-linear}, the first term is written as follows
$$
\int_0^\tau \int_{\R^3} 
(\varrho - 1)\mathbf{u}^E \cdot \partial_t \nabla \Phi_{\eps,\delta}\; dxdt 
$$
$$
= 
\left[
\int_0^\tau \int_{\R^3} 
\gamma(\varrho - 1)\mathbf{u}^E \cdot \left(
-\frac{1}{\varepsilon}\nabla s_{\varepsilon,\delta}
+2\kappa^2\varepsilon\Delta s_{\varepsilon,\delta}
\right)\; dxdt
\right]\to 0 \ \ \ \ \text{as} \ \ \ \ \varepsilon\to0.
$$
We conclude analyzing the last three terms in (\ref{I_1-split}). Owing to Corollary \ref{cor-AW-dec}, we have
$$
\left[
\int_0^\tau \int_{\R^3} 
(\varrho\mathbf{U} - \varrho\mathbf{u}) \otimes \nabla \Phi_{\eps,\delta} : \nabla \mathbf{u}^E \; dxdt
\right.
$$

$$
\left.
+ \int_0^\tau \int_{\R^3} 
(\varrho\mathbf{U} - \varrho\mathbf{u}) \otimes \mathbf{u}^E : \nabla^2 \Phi_{\eps,\delta}\; dxdt
\right.
$$

$$
\left.
+ \frac{1}{2}\int_0^\tau \int_{\R^3} 
(\varrho\mathbf{U} - \varrho\mathbf{u}) \cdot
\nabla |\nabla \Phi_{\eps,\delta}|^2 \; dxdt
\right]
\to 0 \ \ \ \ \text{as} \ \ \ \ \varepsilon\to0.
$$
In particular, the last term writes
\begin{align*}
&\frac{1}{2}\int_0^\tau \int_{\R^3} 
(\varrho\mathbf{U} - \varrho\mathbf{u}) \cdot
\nabla |\nabla \Phi_{\eps,\delta}|^2 \; dxdt\\
&=\frac{1}{2}\int_0^\tau \int_{\R^3} 
\left(\sqrt{\rho}-1\right)\left(\sqrt{\varrho}\mathbf{U} - \sqrt{\varrho}\mathbf{u}\right) \cdot
\nabla |\nabla \Phi_{\eps,\delta}|^2+\left(\sqrt{\varrho}\mathbf{U} - \sqrt{\varrho}\mathbf{u}\right)\cdot \nabla |\nabla \Phi_{\eps,\delta}|^2 \; dxdt.   
\end{align*}
The former converges to $0$ in view of the uniform bounds \eqref{eq: decay reps} and the latter by virtue of Corollary \ref{cor-AW-dec}.

\subsection{Korteweg terms} \label{kt}

The Korteweg terms amount to $I_2$ and $I_3$ as defined in \eqref{step-3}. For a fixed $\kappa > 0$ and from (\ref{ansatz}) we have for $I_2$ that 
\begin{align*}
I_2=&2\kappa^2 \int_0^{\tau} \int_{\R^3} (\varrho\partial_t \Delta r + \varrho\mathbf{u} \cdot \nabla\Delta r-\partial_t r \Delta r)dxdt\\
&=\int_0^{\tau} \int_{\R^3}
(\varepsilon\varrho\partial_t \Delta s_{\varepsilon,\delta} + \varepsilon(\varrho\mathbf{u}) \cdot \nabla\Delta s_{\varepsilon,\delta}
-\varepsilon^2\partial_t s_{\varepsilon,\delta} \Delta s_{\varepsilon,\delta}
)dxdt \to 0 \ \ \ \ \text{as} \ \ \ \ \varepsilon \to 0
\end{align*}
thanks to the uniform bounds (\ref{bounds-1}), (\ref{bounds-2}) for $\rho-1$ and $\rho u$ respectively as well as the bounds \eqref{Hs-cons} and \eqref{eq:Strichartz final} for $s_{\eps,\delta}, \Phi_{\eps,\delta}$.
Upon using \eqref{ansatz} and $\div\ue=0$, we write $I_3$ as
\begin{align*}
I_3&= -2\kappa^2 \int_0^{\tau} \int_{\R^3} 
\left( 
 \nabla\varrho \cdot \nabla( \varrho \textrm{div}\mathbf{U}) 
 - \frac{1}{2} |\nabla \varrho|^2 \textrm{div}\mathbf{U}
 + \nabla \varrho \otimes \nabla \varrho : \nabla \mathbf{U}
 \right)dxdt.\\
&=-\kappa^2\int_0^{\tau} \int_{\R^3}
|\nabla \varrho|^2\Delta \Phi_{\eps,\delta}dxdt-2\kappa^2\int_0^{\tau} \int_{\R^3} 
\varrho\nabla\varrho \cdot \nabla\Delta \Phi_{\eps,\delta}dxdt\\
&-2\kappa^2\int_0^{\tau} \int_{\R^3}
\nabla \varrho \otimes \nabla \varrho : \nabla (\mathbf{u}^E + \nabla \Phi_{\eps,\delta}))
dxdt\\
&= I_3^{(1)} + I_3^{(2)} + I_3^{(3)}.
\end{align*}
The integral $I_3^1$ is bounded as 
\begin{align*}
    &\left|\kappa^2\int_0^{\tau} \int_{\R^3}
\frac{1}{2} |\nabla \varrho|^2 \Delta \Phi_{\eps,\delta}dxdt\right|
\leq \kappa^2\int_0^{\tau} \int_{\R^3}|\nabla\varrho-\nabla r|^2 \left|\Delta \Phi_{\eps,\delta}\right|+\left|\nabla r\right|^2 \left|\Delta \Phi_{\eps,\delta}\right|dxdt\\
&\leq \kappa^2\|\Delta\Phi_{\eps,\delta}\|_{L^{\infty}([0,T]\times \R^3)} \int_0^{\tau} \int_{\R^3}|\nabla\varrho-\nabla r|^2 dxdt+
\|\nabla s_{\eps,\delta}\|_{L^{2}L^{4}}^2\|\Delta\Phi_{\eps,\delta}\|_{L^{2}L^2}\\
&\leq C\int_0^{\tau}\mathcal{E}(\cdot,t)d t+C\kappa^2 T^{\beta} \eps^{\alpha}
\end{align*}
for some $\alpha,\beta >0$ small from \eqref{eq:Strichartz final}. Further and due to \eqref{Hs-cons}, for any $s>\frac{5}{2}$ it holds
\begin{equation}\label{eq:Linfty-bound}
    \|\Delta\Phi_{\eps,\delta}\|_{L^{\infty}([0,T]\times \R^3)}\leq\|\nabla\Phi_{\eps,\delta}\|_{L^{\infty}H^s}\leq C.
\end{equation}
The term  $I_3^2$ obeys the bound
\begin{align*}
&\left|\int_0^\tau\int_{\R^3}\rho\nabla\rho\nabla\Delta\Phi_{\eps,\delta}dxdt\right|=\left|\int_0^{\tau}\int_{\R^3}(\rho-1)\nabla\rho\nabla\Delta\Phi_{\eps,\delta}+\nabla (\rho-1)\nabla\Delta\Phi_{\eps,\delta}dx dt\right|\\
&\leq T^{\beta} \|\rho-1\|_{L^{\infty}L^4}\|\nabla \rho\|_{L^{\infty}L^2}\|\nabla\Delta\Phi_{\eps,\delta}\|_{L^{\frac{8}{3}}L^4}+\left|\int_0^\tau\int_{\R^3}\nabla(\rho-1)\nabla\Delta\Phi_{\eps,\delta}dx dt\right|\\
&\leq CT^{\beta}\eps^{\alpha}+\left|\int_0^\tau\int_{\R^3}\nabla(\rho-1)\nabla\Delta\Phi_{\eps,\delta}dx dt\right|.
\end{align*}
Note that as $\nabla(\rho-1)\in L^{\infty}(0,T;H^1(\R^3))$ is uniformly bounded and $\rho-1\rightarrow 0$ strongly in $L^{\infty}(0,T;H^s(\R^3))$ for all $0\leq s<1$ from \eqref{bounds-1}. It follows that $\nabla(\rho-1)\rightharpoonup^{\ast} 0$ in $L^{\infty}(0,T;L^2(\R^3))$ as $\eps\rightarrow 0$. Hence,
\begin{equation*}
\int_0^\tau\int_{\R^3}\nabla(\rho-1)\nabla\Delta\Phi_{\eps,\delta}dx dt\rightarrow 0 \quad \text{as} \, \eps\rightarrow 0.
\end{equation*}
Finally, to bound $I_3^{(3)}$, one observes that
\begin{align*}
    |I_3^{(3)}|&\leq\kappa^2\int_0^{\tau} \int_{\R^3}\left|\nabla\rho\right|^2\left|\nabla \mathbf{U}\right|dx dt\\
    &\leq 2\kappa^2\left\|\nabla \mathbf{U}\right\|_{L^{\infty}([0,T]\times \R^3)} \int_0^{\tau} \int_{\R^3}\left|\nabla \rho-\nabla r\right|^2+2\kappa^2\int_0^{\tau} \int_{\R^3}\left|\nabla r\right|^2\left|\nabla\mathbf{U}\right|dx dt.
\end{align*}
Using that $\mathbf{U}=\ue+\nabla\Phi_{\eps,\delta}$ is smooth and arguing as for $I_{3}^{(1)}+I_{3}^{(2)}$ we conclude that
\begin{equation*}
    |I_{3}^{(3)}|\leq C \int_0^{\tau}\mathcal{E}(\cdot,t)d t+C\kappa T^{\beta} \eps^{\alpha}
\end{equation*}
for some $\alpha,\beta>0$.

\subsection{Viscous terms} \label{vt}

We have
$$
I_4 = 2\nu \int_0^\tau \int_{\R^3} \sqrt{\varrho}
\mathcal{S}(\mathbf{u})
: \nabla \mathbf{u}^E  dxdt +
2\nu\int_0^\tau \int_{\R^3} \sqrt{\varrho}
\mathcal{S}(\mathbf{u})
: \nabla \nabla \Phi_{\eps,\delta}\; dxdt = I_4^{(1)} + I_4^{(2)}.
$$
Now,
$$
I_4^{(1)} = 2\nu \int_0^\tau \int_{\R^3} (\sqrt{\varrho} - 1)
\mathcal{S}(\mathbf{u})
: \nabla \mathbf{u}^E  dxdt
+
2\nu \int_0^\tau \int_{\R^3}
\mathcal{S}(\mathbf{u})
: \nabla \mathbf{u}^E  dxdt.
$$

$$
\leq C \sqrt{\nu}
\left[
\int_0^\tau \int_{\R^3} |\sqrt{\varrho} - 1|^2
dx dt\right]^{1/2}
\sqrt{\nu}
\left[
\int_0^\tau \int_{\R^3} |\mathcal{S}(\mathbf{u})|^2
dxdt\right]^{1/2}
$$

$$
+ C \sqrt{\nu}
\left[
\int_0^\tau \int_{\R^3} |\mathcal{S}(\mathbf{u})|^2
dxdt\right]^{1/2}
\sqrt{\nu}
$$

$$
\leq \frac{C}{2} \nu
\int_0^\tau \int_{\R^3} |\sqrt{\varrho} - 1|^2
dx dt
+
\frac{C}{2} \nu
\int_0^\tau \int_{\R^3} |\mathcal{S}(\mathbf{u})|^2
dxdt
$$

$$
+ \frac{C}{2} \nu
\int_0^\tau \int_{\R^3} |\mathcal{S}(\mathbf{u})|^2
dxdt
+
C(\nu).
$$
The first term goes to zero as $\varepsilon \to 0$ thanks to (\ref{bounds-1}). The viscous terms containing $\mathcal{S}(\mathbf{u})$ can be absorbed on the left-hand side of (\ref{step-3}) and $C(\nu) \to 0$ as $\nu \to 0$.
For the second term, we have
$$
I_4^{(2)} = 2\nu \int_0^\tau \int_{\R^3} (\sqrt{\varrho} - 1)
\mathcal{S}(\mathbf{u})
: \nabla \nabla \Phi_{\eps,\delta}\;  dxdt
+
2\nu \int_0^\tau \int_{\R^3}
\mathcal{S}(\mathbf{u})
: \nabla \nabla \Phi_{\eps,\delta} \; dxdt
$$

$$
\leq
C\nu
\| \sqrt{\varrho} - 1 \|_{L^4_t L^6_x}
\| \mathcal{S}(\mathbf{u}) \|_{L^2_t L^2_x}
\| \nabla \nabla \Phi_{\eps,\delta} \|_{L^4_t L^3_x}
$$

$$
+
C
\sqrt{\nu}
\left[
\int_0^\tau \int_{\R^3} |\mathcal{S}(\mathbf{u})|^2
dxdt\right]^{1/2}
\sqrt{\nu}
\left[
\int_0^\tau \int_{\R^3} |\nabla \nabla \Phi_{\eps,\delta}|^2
dxdt\right]^{1/2}
$$

$$
\leq
C\eps^{\alpha}\sqrt{\nu}+
C
\sqrt{\nu}
\left[
\int_0^\tau \int_{\R^3} |\mathcal{S}(\mathbf{u})|^2
dxdt\right]^{1/2}
\sqrt{\nu}
\left[
\int_0^\tau \int_{\R^3} |\nabla \nabla \Phi_{\eps,\delta}|^2
dxdt\right]^{1/2}, \ \ \ \ (0 \leq \alpha \leq 1/12).
$$
Here we used \eqref{eq: decay reps}, 
\eqref{eq:Strichartz final} and the second term on the right-hand side can be handled similarly as in $I_4^{(1)}$.

\subsection{Pressure terms} \label{pt}

We recall
$$
I_5 = - \frac{1}{\varepsilon^2} \int_0^\tau\!\! \int_{\R^3}  \left[  \partial_t  (H'(r))(\varrho - r) 
+ \varrho\mathbf{u} \cdot \nabla (H'(r)) +p(\varrho) \mbox{div} \mathbf{U} \right]  dxdt.
$$
We first consider the second term in $I_5$. Using $r=1+\varepsilon s_{\eps,\delta}$, we have
$$
- \frac{1}{\varepsilon^2} \int_0^\tau\!\! \int_{\R^3}
\varrho\mathbf{u} \cdot \nabla (H'(r))
dxdt
= 
- \frac{1}{\varepsilon} \int_0^\tau\!\! \int_{\R^3}
\varrho\mathbf{u} \cdot H''(r) \nabla s_{\eps,\delta}
dxdt
$$

$$
= -\int_0^\tau\!\! \int_{\R^3}
\varrho\mathbf{u} \cdot \nabla s_{\eps,\delta}
\frac{H''(1+\varepsilon s_{\eps,\delta}) - H''(1)}{\varepsilon}
dxdt
- \frac{1}{\varepsilon} \int_0^\tau\!\! \int_{\R^3}
\varrho\mathbf{u} \cdot H''(1) \nabla s_{\eps,\delta}
dxdt,
$$
where $H''(1) = p'(1) = \gamma$. Realizing that 
$$
\left|
\frac{H''(1+\varepsilon s_{\eps,\delta}) - H''(1)}{\varepsilon}
\right|
\leq C |s_{\eps,\delta}|,
$$
we have
$$
\left[
\int_0^\tau\!\! \int_{\R^3}
\varrho\mathbf{u} \cdot \nabla s_{\eps,\delta}
\frac{H''(1+\varepsilon s_{\eps,\delta}) - H''(1)}{\varepsilon}
dxdt
\right] \to 0 \ \ \ \ \text{as} \ \ \ \ \varepsilon\to0,
$$
from Corollary \ref{cor-AW-dec}. For the other term, using the acoustic equation \eqref{eq:AW-linear},
$$
- \frac{\gamma}{\varepsilon} \int_0^\tau\!\! \int_{\R^3}
\varrho\mathbf{u} \cdot 
\nabla s_{\eps,\delta}
dxdt = 
\int_0^\tau \int_{\R^3} 
\varrho\mathbf{u} \cdot \partial_t \nabla \Phi_{\eps,\delta}\; dxdt - 
2\gamma\eps\kappa^2\int_0^\tau \int_{\R^3} 
\varrho\mathbf{u} \cdot\nabla\Delta s_{\eps,\delta} \; dxdt
$$
where, without loss of generality, we assumed $H''(1) = 1$.
The first term cancels with its counterpart $I^{I}$ in \eqref{I-II} while the second converges to $0$. Now, 

$$
\frac{1}{\varepsilon^2} \int_0^\tau\!\! \int_{\R^3}    \partial_t  (H'(r))(r - \varrho) 
- p(\varrho) \mbox{div} \mathbf{U}   dxdt
$$
$$
= \int_0^\tau\!\! \int_{\R^3}
\frac{1-\varrho}{\varepsilon}H''(r)\partial_t s_{\eps,\delta}\; dxdt
+ \int_0^\tau\!\! \int_{\R^3}
s_{\eps,\delta} H''(r) \partial_t s_{\eps,\delta}\; dxdt
-\int_0^\tau\!\! \int_{\R^3} \frac{p(\varrho)}{\eps^2}\Delta \Phi_{\eps,\delta}dxdt
$$
$$= \int_0^\tau\!\! \int_{\R^3}\frac{1-\varrho}{\varepsilon} H''(1) \partial_t s_{\eps,\delta}+ \int_0^\tau\!\! \int_{\R^3}\frac{1-\varrho}{\varepsilon}(H''(r)-H''(1))\partial_t s_{\eps,\delta}dxdt
$$
$$+ \int_0^\tau\!\! \int_{\R^3}
s_{\eps,\delta} H''(r) \partial_t s_{\eps,\delta}\; dxdt
-\int_0^\tau\!\! \int_{\R^3} \frac{p(\varrho)}{\eps^2}\Delta \Phi_{\eps,\delta} 
$$
$$= \int_0^\tau\!\! \int_{\R^3} 
s_{\eps,\delta} H''(r) \partial_t s_{\eps,\delta}\; dxdt +\int_0^\tau\!\! \int_{\R^3}\frac{1-\varrho}{\varepsilon}(H''(r)-H''(1))\Delta \Phi_{\eps,\delta}dxdt
$$
$$-\int_0^\tau\!\! \int_{\R^3}\frac{p(\varrho)-p'(1)(\varrho-1)-p(1)}{\eps^{2}}\Delta \Phi_{\eps,\delta}dxdt.
$$

Consequently, for these remaining terms we have
$$
\left[
- \int_0^\tau\!\! \int_{\R^3}
\frac{p(\varrho)-p'(1)(\varrho-1)-p(1)}{\varepsilon^2}
\Delta \Phi_{\eps,\delta} \; dxdt
\right.
$$

$$
\left.
-
\int_0^\tau\!\! \int_{\R^3}
\frac{1-\varrho}{\varepsilon}\frac{(H''(r)-H''(1))}{\varepsilon}\Delta\Phi_{\eps,\delta}\; dxdt
\right] \to 0 \ \ \ \ \text{as} \ \ \ \ \varepsilon\to0,
$$
by combining \eqref{bounds-1} for $\gamma\geq 2$ and \eqref{eq:Orlicz} for $1<\gamma<2$ and Corollary \ref{cor-AW-dec}.
Finally,
$$
\int_0^\tau\!\! \int_{\R^3}
s_{\eps,\delta} H''(r) \partial_t s_{\eps,\delta}\; dxdt
$$

$$
= \int_0^\tau\!\! \int_{\R^3}
s_{\eps,\delta} H''(1) \partial_t s_{\eps,\delta}\; dxdt
+
\int_0^\tau\!\! \int_{\R^3}
s_{\eps,\delta} (H''(r)-H''(1)) \partial_t s_{\eps,\delta}\; dxdt
$$
where, using similar arguments as above,
$$
\left[
\int_0^\tau\!\! \int_{\R^3}
s_{\eps,\delta} (H''(r)-H''(1)) \partial_t s_{\eps,\delta}\; dxdt
\right] \to 0 \ \ \ \ \text{as} \ \ \ \ \varepsilon\to0
$$
and
$$
\int_0^\tau\!\! \int_{\R^3}
s_{\eps,\delta} H''(1) \partial_t s_{\eps,\delta}\; dxdt
=
\frac{1}{2} \int_{\R^3} 
|s_{\eps,\delta}|^2
dx|_0^\tau.
$$
Back to $I^{II}$ in \eqref{I-II}, we have
$$
\frac{1}{2} \int_0^\tau \int_{\R^3} \varrho\partial_t|\nabla \Phi_{\eps,\delta}|^2 \; dxdt
$$
$$
= \frac{1}{2} \int_0^\tau \int_{\R^3} (\varrho - 1)\partial_t|\nabla \Phi_{\eps,\delta}|^2 \; dxdt
+
\frac{1}{2} \int_{\R^3} |\nabla \Phi_{\eps,\delta}|^2 \; dx |_0^\tau.
$$
By acoustic energy conservation, we have
$$
\frac{1}{2} \int_{\R^3} 
|s_{\eps,\delta}|^2
dx|_0^\tau
+
\frac{1}{2} \int_{\R^3} |\nabla \Phi_{\eps,\delta}|^2 \; dx |_0^\tau =0
$$
and
$$
\frac{1}{2} \int_0^\tau \int_{\R^3} (\varrho - 1)\partial_t|\nabla \Phi_{\eps,\delta}|^2 \; dxdt
=
$$
$$
\frac{1}{2} \int_0^\tau \int_{\R^3} (\varrho-1)\nabla \Phi_{\eps,\delta} \cdot\left(-\frac{1}{\eps}\nabla s_{\eps,\delta}+2\eps\kappa^2\nabla\Delta s_{\eps,\delta}\right)\; dxdt \to 0 \ \ \ \ \text{as} \ \ \ \ \varepsilon\to0
$$
from \eqref{bounds-1} for $\gamma\geq 2$ and \eqref{eq:Orlicz} for $1 <\gamma<2$ as well as Corollary \ref{cor-AW-dec}.
\subsection{Proof of Theorem \ref{main}}\label{proof-closing}

From (\ref{step-3}) and all the estimates above, we end up with
\begin{equation}
\label{step-fnal-1}
\mathcal{E}(\tau)
\leq \mathcal{E}(0)+C \int_0^{\tau}\mathcal{E}(\cdot,t)d t
+ \chi(\nu,\eps)
\end{equation}
where $\chi(\nu,\eps) \to 0$ as $(\nu,\eps) \to 0$ (for fixed $\delta > 0$). In virtue of the integral form of the Gronwall's inequality, we have
\begin{equation}
\label{step-fnal-2}
\mathcal{E}(\tau)
\leq C(T)\mathcal{E}(0)+\chi(\nu,\eps).
\end{equation}
Consequently, sending $\eps \to 0$ first, and then $\delta \to 0$, thanks to \eqref{eq:limitData-2}, \eqref{eq:limitData} and \eqref{conv-id}, we have
\begin{equation}
\label{step-fnal-3}
\lim_{\delta \to 0} \lim_{\eps \to 0} \mathcal{E}(\tau) = 0 
\end{equation}
uniformly in $\tau \in (0,T)$.
Now, considering $\mathbf{U} = \mathbf{u}^E + \nabla \Phi_{\eps,\delta}$, for any compact set $K \subset  \mathbb{R}^3$ we have
\begin{equation}
\label{step-fnal-4}
\|\sqrt{\varrho}(\mathbf{u}-\mathbf{u}^E) \|_{L_t^{2}L_x^2(K)}
\leq
\|\sqrt{\varrho}(\mathbf{u}-\mathbf{U}) \|_{L_t^{2}L_x^2(K)}
+
\|\sqrt{\varrho}\nabla \Phi_{\eps,\delta}\|_{L_t^{p}L_x^q(K)}
\end{equation}
for $p,q > 2$.
The first quantity on the right-hand side of (\ref{step-fnal-4}) goes to zero thanks to (\ref{step-fnal-3}), while
\begin{equation} \label{last-conv}
\lim_{\delta \to 0} \lim_{\eps \to 0}
\|\sqrt{\varrho}\nabla \Phi_{\eps,\delta}\|_{L_t^{p}L_x^q(K)}=0
\end{equation}
thanks to (\ref{eq: decay reps}),  (\ref{eq:Strichartz final}) and (\ref{Hs-cons}) provided that $(p,q)$ is a Schr\"odinger admissible pair. In fact, we have
$$
\|\sqrt{\varrho}\nabla \Phi_{\eps,\delta}\|_{L_t^{p}L_x^q(K)}
\leq \|(\sqrt{\varrho} - 1)\nabla\Phi_{\eps,\delta}\|_{L_t^{p}L_x^q(K)}
+
\|\nabla\Phi_{\eps,\delta}\|_{L_t^{p}L_x^q(K)}
\to 0 \ \ \ \ \text{as} \ \ \ \ \varepsilon\to0
$$
for a fixed $\delta > 0$. Consequently, (\ref{last-conv}) holds. Theorem \ref{main} is proved.

\section{Comments on the 2D case}
\label{2D-c}
In this section, we discuss an extension of our main result to the problem posed on $\R^2$. To that end, we only highlight the key steps and modifications needed for its proof. The relevant differences compared to $d=3$ are given by
\begin{itemize}
    \item the properties of the target system, namely the $2D$ Euler equations, 
    \item while $\varrho_\epsilon-1$ enjoys slightly better integrability properties due to Sobolev embedding, the dispersion of the acoustic waves is weaker in $\R^2$ compared to $\R^3$.
\end{itemize}
To start and concerning the $2D$ incompressible Euler equations, one has global existence of strong solutions as singularity formation in finite time is ruled out by the Beale-Kato-Majda criterion and the conservation of the $L^{\infty}$-norm of the vorticity, see e.g. \cite{bertozzi}.
\begin{theorem}\label{2dEuler-ex}
Given $\mathbf{u}_0^E\in W^{3,2}(\R^2)$ with ${\rm div}\mathbf{u}_0 = 0$, there exists a unique solution
$$\mathbf{u}^E\in C^k([0,\infty),W^{3-k,2}(\R^2;\R^2)), \,\Pi\in C^k([0,\infty),W^{3-k,2}(\R^2)),\, k=0,1,2,3$$
to the initial value problem
(\ref{EULER}) - (\ref{EULER-ic})
such that for all $0<T<\infty$ it holds
\begin{equation}\label{2Deeuler}
\|\mathbf{u}^E\|_{W^{k,\infty}(0,T;W^{3-k,2}(\R^3;\R^3))} + \|\Pi\|_{W^{k,\infty}(0,T;W^{3-k,2}(\R^3))}\leq c(T,\|\mathbf{u}_0^E\|_{W^{3,2}(\R^3)}).
\end{equation}
\end{theorem}
For the primitive system \eqref{mom} posed on $[0,\infty)\times \R^2$, we postulate existence of global finite energy weak solutions according to Definition \ref{weak-def} adapted to $d=2$. Similarly to the case $d=3$, we expect that global existence can be inferred in the spirit of Theorem \ref{main-ex} by relying on \cite{PaoloStefano, LarsPaoloStefano}.

\begin{theorem} \label{main-ex-2D}
Given initial data $(\varrho^0, \mathbf{u}^0)$ of finite energy and BD-entropy, then, there exists at least a global weak solution $(\varrho, \mathbf{u}, \mathcal{T})$ of \eqref{mom} posed on $[0,T)\times \R^2$ in the sense of Definition \ref{weak-def}. 
\end{theorem}
The constructed finite energy weak solutions satisfy the following uniform bounds. 

\begin{lemma}\label{lem:initial data-2D}
Let $(\varrho_\epsilon, \mathbf{u}_\eps^0)$ be a FEWS to \eqref{mom} posed on $[0,T)\times \R^2$ with initial data $(\varrho_\eps^0,\mathbf{u}_\eps^0)$ satisfying  $$E(\varrho_\eps^0,\mathbf{u}_\eps^0)\leq C, \quad   B(\varrho_\eps^0,\mathbf{u}_\eps^0)\leq C,$$
for some $C>0$ independent of $\epsilon>0$. Then, the following hold:
\begin{enumerate}[(i)]
    \item $\varrho_\varepsilon^0-1\in H^1(\R^2)$ and $\varrho_\epsilon-1\in L^{\infty}(0,T;H^1(\R^2))$ uniformly bounded and 
    \begin{equation*}
          \|\varrho_\varepsilon^0-1\|_{L^2(\R^2)}\leq C\eps, \qquad  \|\varrho_\varepsilon-1\|_{L^{\infty}(0,T;L^2(\R^2))}\leq C\eps,
    \end{equation*}
    \item $\sqrt{\varrho_\varepsilon^0}-1\in L^2(\R^2)$, $\sqrt{\varrho_\varepsilon}-1\in L^{\infty}(0,T;L^2(\R^2))$  uniformly bounded with
         \begin{equation*}
          \|\sqrt{\varrho_\varepsilon^0}-1\|_{L^2(\R^2)}\leq C\eps, \qquad  \|\sqrt{\varrho_\varepsilon}-1\|_{L^{\infty}(0,T;L^2(\R^2))}\leq C\eps,
    \end{equation*}
        \item $\sqrt{\varrho_\varepsilon^0}\mathbf{u}_\eps^0\in L^2(\R^2)$ uniformly bounded,
    \item the momentum satisfies $\varrho_\varepsilon^0\mathbf{u}_\eps^0\in L^2(\R^2)+L^{p}(\R^2)$, $\varrho_\varepsilon^0\mathbf{u}_\eps^0\in L^{\infty}(0,T;L^2(\R^2)+L^{p}(\R^2))$ uniformly bounded for all $1\leq p<2$,
    \item the density fluctuations $\seps^0$ satisfy  $\seps^0\in L^2(\R^2)$, $\seps\in L^{\infty}(0,T;L^2(\R^2))$ uniformly bounded as well as $\eps\seps^0\in H^1(\R^2)$, $\eps\seps\in L^{\infty}(0,T;H^1(\R^2))$ uniformly bounded.
\end{enumerate}
\end{lemma}

\begin{proof}
The only modification in the proof compared to one of Lemma \ref{lem:initial data} consists in the fact that the argument yielding $\varrho_\eps-1\in L^6(\R^3)$ from the bound $\nabla\varrho_\epsilon\in L^2(\R^3)$ does not hold for $d=2$. Instead, we rely on the fact that if $f$ is a measurable function such that $\nabla f\in L^2(\R^2)$ and $\supp(f)$ is of finite Lebesgue measure then
\begin{equation}\label{ineq:support}
    \|f\|_{L^p(\R^2)}\leq \|\nabla f\|_{L^2(\R^2)}\mathcal{L}^2\left(\supp(f)\right)^{\frac{1}{p}}.
\end{equation}
For a proof of \eqref{ineq:support} see for instance \cite[Inequality (3.10)]{BrezisLieb} and also \cite[Proof of Lemma 3.1]{lars2}. From the analogue of \eqref{eq:Orlicz} for $d=2$ and the Chebychev inequality, we infer that 
\begin{equation*}
    \mathcal{L}^2(\{|\reps-1|> \frac12\})\leq \frac{1}{2^{\gamma}} \int_{\R^2}\left|\reps(t)-1\right|^{\gamma}\mathrm{1}_{\{|\reps-1|> c\}}\dd x\leq C\eps^{2},
\end{equation*}
where $\mathcal{L}^2$ denotes the Lebesgue measure. Consider a smooth cut-off $\chi\in C_c^{\infty}(\R)$ such that $\mathbf{1}_{[3/4,5/4]}(r)\leq \chi(r)\leq \mathbf{1}_{[1/2,3/2]}(r)$. Applying the inequality \eqref{ineq:support} to $(\varrho_\eps-1)(1-\chi(\reps))$ yields
\begin{multline*}
        \left\|(\varrho_\eps-1)(1-\chi(\reps))\right\|_{L^{\infty}(0,\infty;L^2(\R^2))}
        \leq \|\nabla\varrho_\eps\|_{L^{\infty}(0,\infty;L^2(\R^2))}\mathcal{L}^2\left(\supp(1-\chi(\reps))\right)^{\frac12}\\
        \leq C \eps\|\nabla\varrho_\eps\|_{L^{\infty}(0,\infty;L^2(\R^2))}\leq C\eps.
\end{multline*}
We obtain that 
\begin{equation*}
    \|\varrho_\epsilon-1\|_{L^{\infty}(0,\infty;L^2(\R^2))}\leq C\eps, \quad   \|\nabla(\varrho_\epsilon-1)\|_{L^{\infty}(0,\infty;L^2(\R^2))}\leq C.
\end{equation*}
The proof for the respective bounds of $\varrho_\epsilon^0-1$ follows verbatim.
\end{proof}

Note that opposite to Lemma \ref{lem:initial data}, none of these bounds depends on whether $\gamma<2$ or $\gamma\geq 2$ for $d=2$. 
In particular, by interpolation it follows that 
\begin{equation*}
    \|\rho_\eps\|_{L^{\infty}(0,T;H^s(\R^2))}\leq C\eps^{1-s}
\end{equation*}
and for all $2\leq q<\infty$ there exists $C>0$ such that 
\begin{equation}\label{eq: decay reps2D}
     \begin{aligned}
         \|\sqrt{\varrho_\varepsilon^0}-1\|_{L^q}&\leq \|\varrho_\varepsilon^0-1 \|_{L_x^q}\leq C \eps^{\frac{2}{q}}\\
         \|\sqrt{\varrho_\varepsilon}-1\|_{L_t^{\infty}L_x^q}&\leq \|\varrho_\varepsilon-1 \|_{L_t^{\infty}L_x^q}\leq C \eps^{\frac{2}{q}}
     \end{aligned}
\end{equation}
In addition, the analogue of Lemma \ref{lemm-BD} remains valid for $d=2$.

Next, we provide the suitable decay of acoustic waves in space-time norms. The dispersion relation \eqref{F-mult} is non-homogeneous and mimics wave-like behavior for low and Schr\"odinger like behavior for high frequencies. While it exhibits a regularizing effect for low frequencies for $d>2$ providing decay of order $\eps^\delta$  at the expense of a loss or regularity of order $\delta$ with $\delta>0$ arbitrarily small, no such regularizing effect occurs for $d=2$, see \cite[Section 3]{lars3}. However, separating frequencies above and below the threshold $\frac{1}{\eps}$ and relying on the wave-like estimate for low frequencies combined with an interpolation argument, the following Strichartz estimates are shown to hold in \cite[Proposition 3.8]{lars3}.
\begin{definition}\label{def:betagamma}
The exponents $(q, r)$ are said to be \emph{$\theta$-admissible} if $2\leq q, r\leq\infty$, $(q, r, \theta)\neq(2, \infty, 1)$ and 
\begin{equation*}
\frac1q+\frac{\theta}{r}=\frac{\theta}{2}.
\end{equation*}
We say that a pair is Schr{\"o}dinger or wave admissible if $\theta=\frac{d}{2}$ or $\theta=\frac{d-1}{2}$ respectively.  
Further, we denote $\beta=\beta(r):=\frac12-\frac1r$. 
\end{definition}
\begin{proposition}[\cite{lars3}]\label{prop:Strichartz_hom}
Let $\eps>0$ and $\theta\in[0,1)$. Then, for any $\frac{2-\theta}{2}$-admissible pair $(q,r)$ and $s_0=3\beta(r)\theta$, it holds
\begin{equation}\label{eq:Strichartz-h}
    \|\eith f\|_{L^q(0,T;L^r(\R^2))}\leq C \eps^{\frac{s}{3}}\|f\|_{\dot{H}^{s}(\R^2)}.
\end{equation}
\end{proposition}
The estimate allows one to infer decay in Strichartz-norms at the cost of arbitrarily small regularity. 

\begin{remark}
Being the symbol $\phie$ non-homogeneous, it does not allow for a separation of scales. Therefore, the $\eps$-dependent estimates cannot be obtained by a simple scaling argument. Further, such estimates also appear in the context of the Gross-Pitaevskii equation. For $\theta=0$, \eqref{prop:Strichartz_hom} yields Schr\"odinger like Strichartz estimates that however do not provide decay in $\eps$. In \cite[Corollary B.1]{BDS10}, the estimate \eqref{eq:Strichartz-h} is proven for $\theta=1$ and $d\geq 2$ and low frequencies in the framework of the (GP)-equation. For high frequencies, a Schr\"odinger type estimate is obtained. In this regard, \eqref{eq:Strichartz-h} can be interpreted as a refinement of \cite[Corollary B.1 ]{BDS10}.
\end{remark}

With these key ingredients at hand and observing that the proof of the relative entropy method can then be adapted in a straight-forward way, we have the following main result for the asymptotic limit for $d=2$. Specifically and analogue to the case $d=3$, we consider ill-prepared initial data

\begin{equation}\label{eq:limitData-2-2D}
    \varrho(0,\cdot)=\varrho_\eps^0=1+\eps \seps^0,\quad  \seps^0\in L^{\infty}(\R^2)\cap H^{1}(\R^2) \,\, \text{uniformly bounded}, \,\, \seps^0\rightarrow s^0 \quad \text{in} \,\, H^1(\R^2),
\end{equation}
    
 \begin{equation}
 \label{eq:limitData-2D}
 \mathbf{u} (0,\cdot)=\mathbf{u}_\eps^0,\quad  \sqrt{\varrho_\eps^0}\mathbf{u}_\eps^0 \in L^{2}(\R^2), \quad   \sqrt{\varrho_\eps^0}(\mathbf{u}_\eps^0-\mathbf{u}^0)\rightarrow 0\quad \text{in} \quad L^2(\R^2).
\end{equation}
as $\varepsilon \to 0$.

\begin{theorem} \label{main2D}
    Suppose $\nu \to 0$ as $\varepsilon \to 0$. Assume the initial data $(\varrho_\varepsilon^0, \mathbf{u}_\varepsilon^0)$ to be of uniformly bounded
    finite energy and BD-entropy, i.e. there exists $C>0$ independent of $\varepsilon, \nu>0$ such that 
    \begin{equation}\label{eq: uniform bound data-2D}
        E(\varrho_\eps^0,\mathbf{u}_\eps^0)\leq C, \quad   B(\varrho_\eps^0,\mathbf{u}_\eps^0)\leq C.
    \end{equation}
   and satisfy \eqref{eq:limitData-2-2D}, \eqref{eq:limitData-2D}. 
   
   Let $(\varrho_\varepsilon, \mathbf{u}_\varepsilon, \mathcal{T}_\varepsilon)$ be a weak solution of (\ref{mom}) and $\mathbf{u}^E$ the unique solution to the initial value problem (\ref{EULER}) - (\ref{EULER-ic}) on $[0,\infty)\times \R^2$  with initial datum $\mathbf{u}_{0}^E = \mathbf{P}(\mathbf{u}^0) \in W^{3,2}(\mathbb{R}^3)$. Then, for all $T>0$ as $\varepsilon \to 0$, it holds
\begin{equation}
\label{main-rho-conv2D}
\|\varrho_\varepsilon-1\|_{L^\infty (0,T;H^s(\R^2))}\leq C\eps^{1-s} 
\end{equation}
for all $0<s<1$ and
\begin{equation}
\label{main-u-conv2D}
\sqrt{\varrho}_\varepsilon\mathbf{u}_\varepsilon \to \mathbf{u}^E
\ \ \text{in} \ \ L^2(0,T;L^2_{\mathrm{loc}}(\mathbb{R}^2))).
\end{equation}
\end{theorem} 

We note that different to Theorem \ref{main}, we may consider arbitrarily large times $T>0$ as solutions to the limit system are global. The convergence in Orlicz space is implied by the convergence in $\ L^\infty (0,T;H^s(\R^2))$ due to the Sobolev embeddings for $d=2$ and any $\gamma>1$.

\section*{Acknowledgments} M. C. has been supported by the Praemium Academiae of \v S. Ne\v
casov\' a. The Institute of Mathematics CAS is supported by RVO:67985840. D.D. gratefully acknowledge the partial support by the GruppoNa\-zio\-na\-le per l'Analisi Matematica, la Probabilit\`a e le loro Applicazioni (GNAMPA) of the Istituto Nazionale di Alta Matematica(INdAM) and by the by PRIN2022-PNRR- Project N. P20225SP98 ``Some mathematical approaches to climate change and its impacts''. L.E.H. is funded by the Deutsche Forschungsgemeinschaft (DFG, German Research Foundation) – Project- ID 258734477 – SFB 1173.

\section*{Data Availability Statement}
We declare that the manuscript has no associated data and hence no data set has been used in the realization of this work.

\section*{Declaration - conflict of interest}
The authors have no competing interests to declare that are relevant to the content of this article.

\bibliographystyle{siam}
\bibliography{references}

\end{document}